\documentclass[11pt,bezier,amstex]{article}  

\usepackage{color}              
\usepackage{graphicx}
\usepackage[]{amsmath}
\usepackage{amssymb}
\usepackage{hyperref}
\usepackage{enumerate, amsfonts, amsthm, bbm}

\definecolor{MyDarkBlue}{rgb}{0,0.08,0.50}
\definecolor{BrickRed}{rgb}{0.65,0.08,0}

\hypersetup{
colorlinks=true,       
    linkcolor=MyDarkBlue,          
    citecolor=BrickRed,        
    filecolor=red,      
    urlcolor=cyan           
}

\newtheorem{Lemma}{Lemma}[section]
\newtheorem{Proposition}[Lemma]{Proposition}

\newtheorem{Theorem}[Lemma]{Theorem}

\newtheorem{Corollary}[Lemma]{Corollary}

\newcommand{\prob}{\mathbb{P}}

\newcommand{\E}{\mathbb{E}}

\newcommand{\II}{\mathcal{I}}

\renewcommand{\SS}{\mathcal{S}}

\newcommand{\expec}{\mathbb{E}}

\DeclareMathOperator*{\argmin}{arg\,min}
\newcommand{\e}{{\mathrm e}}

\newcommand{\R}{{\mathbb R}}

\newcommand{\eqn}[1]{\begin{equation} #1 \end{equation}}
\newcommand{\eqan}[1]{\begin{align} #1 \end{align}}
\newcommand{\lbeq}[1]{\label{#1}}
\newcommand{\refeq}[1]{(\ref{#1})}
\newcommand{\sss}{\scriptscriptstyle}

\newcommand{\cluster}{{\cal C}}

\newcommand {\vep}{\varepsilon}
\newcommand\1{\mathbbm{1}}
\newcommand{\indic}[1]{\1_{\{#1\}}}
\newcommand{\indicwo}[1]{\1_{#1}}

\newcommand{\nn}{\nonumber}

\newcommand{\thetastar}{\theta^*}

\newcommand{\thetahat}{\hat\theta}
\newcommand{\fff}{k}

\newcommand{\Fvartheta}{\Lambda}

\numberwithin{equation}{section}

\newcommand{\betam}{\tilde{\beta}}

\newcommand{\ii}{{\rm i}}

\definecolor{darkgreen}{rgb}{0,.4,0}
\definecolor{darkagenta}{rgb}{.5,0,.5}
\definecolor{darkred}{rgb}{1,0,0}
\definecolor{darkblue}{rgb}{0,0,.4}

\begin{document}

\author{
Elie A\"id\'ekon\thanks{Laboratoire de Probabilit\'es et Mod\`eles Al\'eatoires,
Universit\'e Paris 6, 4, place Jussieu, 75005 Paris, France. E-mail: {\tt elie.aidekon@upmc.fr}}
\and
Remco van der Hofstad\thanks{Department of Mathematics and
Computer Science, Eindhoven University of Technology, P.O.\ Box 513,
5600 MB Eindhoven, The Netherlands. E-mail: {\tt
rhofstad@win.tue.nl}, {\tt j.s.h.v.leeuwaarden@tue.nl}}
\and
Sandra Kliem\thanks{Fakult\"at f\"ur Mathematik, Universit\"at Duisburg-Esssen, Thea-Leymann-Str. 9, 45127 Essen, Germany. E-mail:  {\tt sandra.kliem@uni-due.de}
}
\and
Johan S.H. van Leeuwaarden$^{\dagger}$
}

\title{Large deviations for power-law thinned L\'evy processes}

\maketitle

\begin{abstract}
This paper deals with the large deviations behavior of a stochastic process called thinned  L\'evy process. This process appeared recently as a stochastic-process limit in the context of critical inhomogeneous random graphs \cite{BhaHofLee09b}. The process has a strong negative drift, while we are interested in the rare event of the process being positive at large times. To characterize this rare event, we identify a tilted measure. This presents some challenges inherent to the power-law nature of the thinned L\'evy process. General principles prescribe that the tilt should follow from a variational problem, but in the case of the thinned L\'evy process this involves a Riemann sum that is hard to control. We choose to  approximate the Riemann sum by its limiting integral, derive the first-order correction term, and prove that the tilt that follows from the corresponding approximate variational problem is sufficient to establish the large deviations results.
\end{abstract}

\vspace{0.3in}


\noindent
{\bf MSC2000 subject classification.}
60C05, 05C80, 90B15.



\section{Introduction}
\label{sec-int}
This paper deals with a relatively new stochastic process called thinned  L\'evy process. This process appeared as a stochastic-process limit in the context of critical inhomogeneous random graphs \cite{BhaHofLee09b} and is also a special case of the multiplicative coalescent  \cite{Aldo97,AldLim98}. In its most basic form, the thinned L\'evy processes  $(\SS_t)_{t\geq 0}$ is  defined as
    \eqn{
    \lbeq{SS-def-ref}
    \SS_t=1+\betam t+\sum_{i=2}^{\infty} c_i [\II_i(t)-c_i t].
    }
Here $ \II_i(t)=\indic{T_i\leq t}$ denotes an indicator process with $T_i$ an exponential random variable with mean $\expec[T_i]=i^{\alpha}$. All $T_i$ are assumed independent. Furthermore, $\betam\in\mathbb{R}$ and we define the coefficients $c_i=i^{-\alpha}$ with $\alpha\in(\frac13,\frac12)$.

 Let us first explain why we have dubbed in \cite{BhaHofLee09b}  this process thinned L\'evy process. Upon replacing $\II_i(t)$ by a Poisson process with rate $c_i$ the process $(\SS_t)_{t\geq 0}$ becomes a spectrally positive L\'evy process, consisting of infinitely many independent Poisson sources and linear drifts. Compared to the Poisson process, the indicator process $\II_i(t)$ only counts the first event, and in that sense thins the  L\'evy process. Note that the sums $\sum_{i=2}^{\infty} c_i \II_i(t)$ and $\sum_{i=2}^{\infty} c_i ^2 t$ cannot be treated separately due to the assumption $\alpha\in(\frac13,\frac12)$. In fact, both sums are connected with the Riemann zeta function defined as $\zeta(s)=\sum_{n=1}^\infty n^{-s}$ for ${\rm Re}(s)>1$,  and for all other values $s\neq 1$ defined via the Riemann functional equation. For the purpose of this paper, it is convenient to work with the analytic continuation
 	\eqn{
 	\lbeq{rz}
 	\zeta(s)=\lim_{N\to\infty}\Big\{\sum_{n=1}^N n^{-s}-\frac{N^{1-s}}{1-s}-\frac12 N^{-s}\Big\},
	\quad {\rm Re}(s)>-1, s\neq 1,
 	}
 which follows from Euler-Maclaurin summation \cite[p.~333]{Har49}, and which gives meaning to the identity $\sum_{i=2}^{\infty} c_i ^2=\zeta(2\alpha)$. Also, using $\expec \indic{T_i\leq t}=1-\e^{-c_i t}$,
     \eqn{
    \lbeq{SS-mean}
    \expec[\SS_t]=1+\betam t+\sum_{i=2}^{\infty} c_i [1-\e^{-c_i t}-c_i t].
    }
Sums of the type as in \eqref{SS-mean} will appear frequently in this paper, and using the results developed in Section \ref{sec-main-term} to replace sums by integrals, it follows that
      	\eqn{
    	\lbeq{SS-mean-ass}
    	\expec[\SS_u]
	\sim u^{-1} \sum_{i \geq 2} g(i/u^{1/\alpha})\sim
	u^{1/\alpha-1}\int_{0}^\infty x^{-\alpha}[1-\e^{-x^{-\alpha}}-x^{-\alpha}]dx,
    	}
where $g(x)=x^{-\alpha}[1-\e^{-x^{-\alpha}}-x^{-\alpha}]$, and the integral expression is finite and negative, so that $\expec[\SS_u]$ decays faster than $u$ for $u\uparrow \infty$ since $\alpha\in(\frac13,\frac12)$. Asymptotics as in \refeq{SS-mean-ass} are made precise in Lemma \ref{lem-taylor}, to which we refer the reader for more details.

The precise power-law form imposed by the assumption $\alpha\in(\frac13,\frac12)$ is essential for our study, not only because it determines the above behavior over time of the mean, but also because this interval for the powers $\alpha$ is intimately related with critical behavior in certain power-law random graphs, as explained next.

In \cite{BhaHofLee09b} scaling limits were obtained for the sizes of the largest components at criticality for rank-1 inhomogeneous random graphs with power-law degrees with power-law exponent $\tau\in (3,4)$ of which we now describe one specific example known as the Poissonian random graph or Norros-Reittu model \cite{NorRei06}. To define the model, we consider the vertex set $[n]:=\{1,2,\ldots, n\}$ and suppose vertex $i$ is assigned a weight $w_i$. Attach an edge between vertices $i$ and $j$ with probability
    \eqn{
    \lbeq{pij-NR}
    p_{ij} = 1-{\rm exp}\Big(-\frac{w_i w_j}{ \sum_{\ell\in[n]} w_\ell}\Big),
    }
Different edges are independent. In this model, the average degree
of vertex $i$ is close to $w_i$, thus incorporating inhomogeneity
in the model. There are many adaptations of this model,
for which equivalent results hold. Indeed,
the model considered here is a special case of the
so-called \emph{rank-1 inhomogeneous random graph}
introduced in great generality in \cite{BolJanRio07}. It is
asymptotically equivalent with many related models (see
\cite{ChuLu02a,ChuLu02b,ChuLu03, ChuLu06c, ChuLu06,BriDeiMar-Lof05}).
See \cite{Jans08a} 
for conditions under which random graphs
are \emph{{asymptotically} equivalent}, meaning that all events
have asymptotically equal probabilities.

Let the weight be defined by
    \eqn{
    \lbeq{choicewi}
    w_i = [1-F]^{-1}(i/n),
    }
where $F$ is a distribution function on $[0,\infty)$
for which we assume that there exists
a $\tau\in (3,4)$ and $0<c_{\sss F}<\infty$ such that $\lim_{x\rightarrow \infty}x^{\tau-1}[1-F(x)]= c_{\sss F}$
with $[1-F]^{-1}(u)=\inf \{ s\colon [1-F](s)\leq u\}$ for $u\in(0,1)$.

For $\tau\in (3,4)$ it was shown that the rescaled sizes of the components converge to hitting times of a thinned L\'evy process.  Let $\cluster(1)$ be the connected component to which the largest-weight node belongs (which is proved to be equal to largest component with high probability). Let $H_1(0)=\inf\{t\geq 0\colon \SS_t=0\}$ denote the first hitting time of $0$ of (a rescaled version of) the process $(\SS_t)_{t\geq 0}$ in \eqref{SS-def-ref} with
	\eqn{
	\alpha:=1/(\tau-1)\in(1/3,1/2).
	}
Below, it will be more convenient to phrase our results in terms of $\tau\in (3,4)$, which we will do from now on.

In \cite[Theorem 2.1]{BhaHofLee09b} it is proved that $|\cluster(1)|n^{-(\tau-2)/(\tau-1)}$ converges in distribution to $H_1(0)$. The critical components are thus of the order $n^{-(\tau-2)/(\tau-1)}$, but  to obtain information beyond the order one needs to investigate $H_1(0)$. In the companion paper \cite{AidHofKliLee13b} we derive the precise asymptotic results for both $\prob(H_1(0)>u)$ and the tail distribution of the largest cluster, for $u\to \infty$. A crucial ingredient of the proofs is the asymptotic behavior of $\prob(\SS_u>0)$, the main topic of the present paper. Indeed, because of the strong downward drift of the process $(\SS_t)_{t\geq 0}$, it seems plausible that, for large $u$, $\prob(H_1(0)>u)\approx \prob(\SS_u>0)$.

We thus study the probability of the rare event $\{S_u>0\}$ for some $u>0$ large. In order to do so, we take the traditional approach to large deviations theory via the so-called change of measure technique,
see e.g.~\cite{DemZei98, Holl00}. In this approach, a tilted measure is identified under which the event $\{S_u>0\}$ has high probability, and the probability of the event under the original measure is estimated in terms of the Radon-Nikodym density relating the two measures. That is, we investigate the measure $\tilde{\prob}$
with Radon-Nikodym derivative $\e^{\vartheta u \SS_{u}}/\expec[\e^{\vartheta u \SS_{u}}]$, for some appropriately chosen $\vartheta$. The choice of $\vartheta$ turns out to be
quite subtle for the thinned L\'evy process. General principles from large deviations theory prescribe that the optimal choice is $\vartheta=\thetahat:=\argmin_{\vartheta}\log \expec[\e^{\vartheta u \SS_u}]$. From \eqref{SS-def-ref} it can be seen that $\log \expec[\e^{\vartheta u \SS_u}]$ is described in terms of an infinite sum that is hard to control. However, this infinite sum is in fact a Riemann sum, which gives rise to the approximation $\log \expec[\e^{\vartheta u \SS_u}]\approx u^{\tau-1}\Fvartheta(\vartheta)$ with $\Fvartheta(\vartheta)$ an integral independent of $u$. Therefore, for large $u$, it should be that
\eqn{
\label{var1}
\thetahat\approx\thetastar:=\arg\min_{\vartheta}[u^{\tau-1}\Fvartheta(\vartheta)]=\arg\min_{\vartheta}\Fvartheta(\vartheta).
}
We could thus apply the tilting with $\thetastar$ instead of $\thetahat$ in the hope to get sharp asymptotic estimates for $\prob(\SS_u>0)$. However, while $\thetastar$ is asymptotically sharp, it turns out to be a too weak approximation of $\thetahat$ for our purposes. We solve this issue by refining the approximative variational problem \eqref{var1} into
\eqn{
\thetahat\approx\thetastar_u:=\arg\min_{\vartheta} [\Fvartheta(\vartheta)
	+\vartheta \vep_u ]
}
with
\eqn{
	\vep_u= \frac{\zeta(\alpha)}{u^{\tau-2}}+ \frac{\betam+1-\zeta(2\alpha)}{u^{\tau-3}},\quad \alpha=\frac{1}{\tau-1}.
	}

The refinement $\vep_u$ that includes the two Riemann zeta functions (defined in \eqref{rz} since $\tau\in(3,4)$) vanishes for $u\to\infty$, and in fact seems only marginal, but it turns out to be crucial in order for the tilting procedure to provide an asymptotically sharp description of the rare event probability $\prob(\SS_u>0)$. This eventually leads to one of the key results of this paper.

\begin{Theorem}[Exact asymptotics tail $\SS_u$]
\label{thm-tail-Su}
There exists $I, D>0$ and $\kappa_{ij}\in {\mathbb R}$
such that, as $u\rightarrow \infty$,
    	\eqn{
    	\lbeq{eqn-tail-C1}
	\prob(\SS_u>0)= \frac{D}{u^{(\tau-1)/2}}
	\e^{-Iu^{\tau-1}+u^{\tau-1}\sum_{i+j \geq 1} \kappa_{ij}u^{-i(\tau-2)-j(\tau-3)}}(1+o(1)).
	}
\end{Theorem}	
Notice that since $\tau\in (3,4)$, the sum over $i,j$ such that $i+j\geq 1$ is in fact \emph{finite}, as we can ignore all terms for which $\tau-1-i(\tau-2)-j(\tau-3)\leq 0$. The asymptotic behavior is dominated by the term $-Iu^{\tau-1}$ with the crucial constant $I$ defined as $I=-\min_{\vartheta\geq 0} \Fvartheta(\vartheta)=-\Fvartheta(\thetastar)>0$. The other constants $D$ and $\kappa_{ij}$ are specified in Sections \ref{sec:main} and \ref{sec:proofld}, and in determining their values it turns out to be crucial to work with the tilting $\thetastar_u$.

In order to derive Theorem \ref{thm-tail-Su}, we shall investigate large deviation properties of $\SS_u$.
The same techniques can be used in order to prove that $\SS_{au}u^{-(\tau-2)}$ approaches a deterministic shape under the conditional distribution given $\SS_u>0$:

\begin{Theorem}[Sample path large deviations]
\label{thm-sample-path-LD}
There exists a function $a\mapsto I_{\sss E}(a)$ on $[0,1]$ such that, for any $\vep>0$ and $a\in [0,1]$,
    	\eqn{
    	\lbeq{eqn-sample-paths}
	\lim_{u\rightarrow \infty}
	\prob\big(\big|\SS_{au}-u^{\tau-2}I_{\sss E}(a)|\leq \vep u^{\tau-2}\mid\SS_u>0)=1.
	}
\end{Theorem}	
See \eqref{def_I_E} for the precise form of $a\mapsto I_{\sss E}(a)$.

\subsection{Discussion}\label{sec-disc}

\paragraph{Large deviations connection.}
By \refeq{SS-mean-ass}, and recalling that we have defined $\alpha=1/(\tau-1)$,
      	\eqn{
    	\lbeq{SS-mean-ass-rep}
    	\expec[\SS_u]
	\sim u^{\tau-2}\int_{0}^\infty x^{-\alpha}[1-\e^{-x^{-\alpha}}-x^{-\alpha}]dx+o(u^{\tau-2})
	\equiv u^{\tau-2}(\mu_{\sss \SS} +o(1)).
    	}
It is not hard to check that $\mu_{\sss \SS}<0$.
Thus, for $u\rightarrow\infty$, the event $\{\SS_u>0\}$ can be thought of as a
\emph{large deviation event}. We next make this connection to large deviation theory more precise.

\paragraph{Classical large deviations.}
We next discuss two connections to classical large deviations. Indeed, when $S_n=\sum_{i=1}^n X_i$,
and $(X_i)_{i=1}^n$ are i.i.d.\ random variables with a finite moment generating function, Cram\'er's Theorem \cite{DemZei98, Holl00} tells us that, for every $a>0$,
	\eqn{
	\label{class-LD-Cramer}
	\prob(S_n-\expec[S_n]\leq -a n)=\e^{-I(a)n}(1+o(1)).
	}
Moreover, by Bahadur-Rao (see e.g.~\cite{DemZei98}),
we have that there exists a constant $A$ such that
	\eqn{
	\label{BR-asymp}
	\prob(S_n-\expec[S_n]\leq -a n)=\frac{A}{n^{1/2}}\e^{-I(a)n}.
	}
Comparing to the main result in Theorem \ref{thm-tail-Su}, we see that a similar result holds with $n$ replaced with $u^{\tau-1}$. This suggests that we can think of Theorem \ref{thm-tail-Su} as describing the classical large deviation result in \eqref{BR-asymp} with $n$ replaced with $u^{\tau-1}$. The only exception is the correction term $u^{\tau-1}\sum_{i+j \geq 1} \kappa_{ij}u^{-i(\tau-2)-j(\tau-3)}$, which is unusual and absent in classical large deviations analysis.

A second connection to large deviations exists with the G\"{a}rtner-Ellis
Theorem \cite{DemZei98, Holl00}. Indeed, in the classical sense, assume again that $S_n=\sum_{i=1}^n X_i$, but now we no longer assume that $(X_i)_{i=1}^n$ are i.i.d.~random variables. Instead, we assume that
	\eqn{
	\label{Ftheta-def}
	\Fvartheta(\vartheta)=\lim_{n\rightarrow \infty}\frac 1n \log\expec[\e^{\vartheta S_n}]
	}
exists. Then,  the G\"{a}rtner-Ellis Theorem tells us that \eqref{class-LD-Cramer} still holds, with
	\eqn{
	I(a)=\sup_{\vartheta} [a\vartheta -\Fvartheta(\vartheta)],
	}
i.e., $a\mapsto I(a)$ is the Legendre transform of $\vartheta \mapsto \Fvartheta(\vartheta)$.
In our setting, we can compute that
	\eqn{
	\label{F-theta-def}
	\Fvartheta(\vartheta)=\lim_{u\rightarrow \infty}\frac 1{u^{\tau-1}} \log\expec[\e^{\vartheta u \SS_u}],
	}
and again $I=\inf_{\vartheta} \Fvartheta(\vartheta)$, which agrees with the G\"{a}rtner-Ellis Theorem
when $a=0$. This explains the philosophy behind the way we have constructed our proof.

\paragraph{Other large deviations events.}
We believe that our methods can be extended to identify the large deviation behavior of other tail events of $\SS_u$, such as $\prob(\SS_u>au^{\tau-2})$ for any $a>\mu_{\sss\SS}$,
where $\mu_{\sss\SS}=\int_{0}^\infty x^{-\alpha}[1-\e^{-x^{-\alpha}}-x^{-\alpha}]dx$ is the asymptotic mean of $u^{-(\tau-2)}\SS_u$ in \refeq{SS-mean-ass-rep}. Alternatively, our methods should extend to events of the form $\prob(\SS_u<au^{\tau-2})$ for any $a<\mu_{\sss \SS}$. Our arguments suggest that such probabilities behave like $\e^{-u^{\tau-1} I_{\sss \SS}(a)(1+o(1))},$ where $I_{\sss \SS}(a)=0$ precisely when $a=\mu_{\sss\SS}$.
In the language of \cite{DemZei98, Holl00}, we expect the random variables $(u^{-(\tau-2)}\SS_u)_{u\geq 0}$ to satisfy a large deviation principle with speed $u^{\tau-1}$ and rate function $I_{\sss \SS}$. The G\"{a}rtner-Ellis Theorem \cite{DemZei98, Holl00} and \eqref{Ftheta-def} then suggests that
	\eqn{
	I_{\sss \SS}(a)=\sup_{\vartheta} [a\vartheta -\Fvartheta(\vartheta)],
	}
where $\vartheta\mapsto \Fvartheta(\vartheta)$ is defined in \eqref{F-theta-def} and computed in \refeq{sum-to-int}.
We do not pursue this further here.

\paragraph{Cluster tails for critical random graphs.}
In \cite{AidHofKliLee13b} we make formal the conjecture that  $\prob(H_1(0)>u)\approx\prob(\SS_u>0)$ for large $u$.
We show that $\prob(H_1(0)>u)$ has the same asymptotic behavior as $\prob(\SS_u>0)$ in \eqref{eqn-tail-C1}, with the same constants except for the constant $D$. Despite the similarity of this result, the proof method in \cite{AidHofKliLee13b} is entirely different. In order to establish the asymptotics for $\prob(H_1(0)>u)$, we establish in  \cite{AidHofKliLee13b}
sample path large deviations, not conditioned on the event $\{\SS_u>0\}$, but on the event $\prob(H_1(0)>u)$.

In particular, in \cite{AidHofKliLee13b}, we establish the following two results. First,
we prove that there exists $A\in(0,D)$
such that
	\eqn{
    	\lbeq{eqn-tail-H1}
	\prob(H_1(0)>u)= \frac{A}{u^{(\tau-1)/2}}
	\e^{-Iu^{\tau-1}+u^{\tau-1}\sum_{i+j \geq 1} \kappa_{ij}u^{-i(\tau-2)-j(\tau-3)}}(1+o(1)).
	}
Equation \eqref{eqn-tail-H1} is much harder than \eqref{eqn-tail-C1} in Theorem \ref{thm-tail-Su}, since
we have to investigate the probability that $\SS_t>0$ for \emph{all} $t\in [0,u]$.
Second, in \cite[Theorem 1.5]{AidHofKliLee13b}, we derive a result related to Theorem \ref{thm-sample-path-LD} saying that
    	\eqn{
    	\lbeq{eqn-sample-paths-strong}
	\lim_{u\rightarrow \infty}
	\prob\big(\sup_{a\in [0,1]}\big|\SS_{au}-u^{\tau-2}I_{\sss E}(a)|\leq \vep u^{\tau-2}\mid H_1(0)>u)=1.
	}
In order to prove \eqref{eqn-sample-paths-strong}, a crucial ingredient is to show that the path cannot deviate much in small time intervals. For this, we need to pay special attention to the fact that time is \emph{continuous}. Indeed, the proof of the extension to \refeq{eqn-sample-paths-strong} consists
of four key steps. In the first, $\big|\SS_{au}-u^{\tau-2}I_{\sss E}(a)|\leq \vep u^{\tau-2}$ with high probability for $a$'s that are close to 0. In the second step, we prove that $\big|\SS_{au}-u^{\tau-2}I_{\sss E}(a)|\leq \vep u^{\tau-2}$ with high probability for a finite, yet growing with $u$, number of values of $a$'s in the interval $[0,1]$ at equal distance that are sufficiently far from the extremeties $a=0$ and $a=1$. In the third step, we show that it is very unlikely that the process $t\mapsto \SS_t$ leaves the tube of width $\vep u^{\tau-2}$ around $u^{\tau-2}I_{\sss E}(a)$ in any of the (small) intervals. In the last and fourth step, we investigate the probability that $\SS_t>0$ for all $t$ close to $u$.
Together, these results suffice to prove \refeq{eqn-sample-paths-strong}.

\paragraph{The case $\tau=4$.}
Although not allowed, it is instructive to substitute $\tau=4$ into \eqref{eqn-tail-C1}. This yields
\eqn{
\lbeq{spec4}
	\prob(\SS_u>0)= \frac{D}{u^{3/2}}
	\e^{-Iu^{3}+\kappa_{01}u^{2}+\kappa_{10}u+\kappa_{11}}(1+o(1)).
	}
This form is reminiscent of results for the
Erd\H{o}s-R\'enyi graph obtained in \cite{HofJanLee10,Pitt01}.
The Erd\H{o}s-R\'enyi graph on the vertex set $[n]:=\{1,\ldots,n\}$ is constructed by including each of the ${n\choose 2}$ possible edges with probability $p$, independently of all other edges. Critical behavior corresponds to $p=(1+\lambda n^{-1/3})/n$,  $\lambda\in\mathbb{R}$ fixed, and letting $n\to\infty$. It is a special case of the rule in  \eqref{pij-NR} when all weights equal $w_i=1+\lambda n^{-1/3}$. Further, for $w_i = [1-F]^{-1}(i/n)$ in \eqref{choicewi}, the same scaling limit for the largest critical clusters holds as for the Erd\H{o}s-R\'enyi random graph when $\expec[W^3]<\infty$, where $W$ has distribution function $F$.

Aldous \cite{Aldo97} showed that the scaling limit describing the critical cluster sizes is a Brownian motion following an asymptotically negative drift of the form $\nu_0+\nu_1t-\nu_2t^2$ with $\nu_2>0$. The size of the largest component, rescaled by $n^{-2/3}$, converges in distribution to some random variable $\gamma_1(\lambda)$. In \cite{HofJanLee10} the excursions of this Brownian motion on a parabola were studied, leading to the result (also derived in \cite{Pitt01} via a different techniques)
  	\eqn{
	\lbeq{pittel1}
    	\prob(\gamma_1(\lambda)>u)=
	\frac{\exp\left(-\frac{1}{8}u(u-2\lambda)^2\right)}{\sqrt{2\pi}u^{3/2}}(1+o(1)),
	\quad u\rightarrow\infty.
   	}
Notice the strong resemblance with \eqref{spec4}. 

\section{Overview of results}\label{sec:main}
In this section we give an overview of the results. Among others, we shall establish Theorem \ref{thm-tail-Su}, announced in the previous section, although this theorem is not the strongest result obtained in this paper. We derive an asymptotic description of the entire density of $\SS_u$ near zero in Proposition \ref{lem-dens-f} from which Theorem \ref{thm-tail-Su} follows, and we extend Theorem \ref{thm-tail-Su} with deriving the optimal trajectory, or sample path large deviations, conditioned on the event $\{\SS_u>0\}$, in Theorem \ref{thm-sample-path-LD}.

Mathematically, establishing these results relies on two main steps. The first step is to consider the variational problem $\min_{\vartheta}\log \expec[\e^{\vartheta u \SS_u}]$ and its minimizer in the asymptotic regime where $u$ is large. In this regime, we can replace the Riemann sum appearing in the expressions for $\log \expec[\e^{\vartheta u \SS_u}]$ by an integral and some first-order correction terms. This gives rise to an asymptotic variational problem that we analyse in great detail using advanced results on bounding sums by integrals and the implicit function theorem. The results are reported in Section \ref{avp}. The second step is to apply the exponential tilting of measure, using the Radon-Nykodym derivative, to establish the properties of the process under the tilted measure. The properties are reported in Section \ref{proppp}. In establishing these properties, it turns out to be sufficient to work with the tilted measure that follows from the solution of the asymptotic variational problem treated in Section \ref{avp}.

\subsection{Asymptotic variational problem}
\label{avp}

We use the notion of \emph{exponential tilting of measure} in order to
give a convenient description of the probability of interest as follows:
    \eqn{
    \lbeq{tilt-form}
    \prob(\SS_{u}>0)
    =\phi(u;\vartheta) \widetilde\expec_{\vartheta}[\e^{-\vartheta u \SS_u}\indic{\SS_{u}>0}],
    }
where $\vartheta$ is chosen later on.
We define the measure $\widetilde \prob_{\vartheta}$ with corresponding expectation
$\widetilde \expec_{\vartheta}$ by the equality, for every event $E$,
    \eqn{
    \lbeq{tildeP-def}
    \widetilde \prob_{\vartheta}(E)=\frac{1}{\phi(u;\vartheta)}\expec_{\vartheta}[\e^{\vartheta u \SS_u}\indicwo{E}],
    }
where the normalizing constant $\phi(u;\vartheta)$ is defined as
    \eqn{
    \lbeq{phi-def}
    \phi(u; \vartheta)=\expec[\e^{\vartheta u \SS_u}].
    }
Choosing a good $\vartheta$ is rather delicate. As discussed around \eqref{F-theta-def},
we would like to choose $\vartheta$ to be the minimizer of
$\vartheta\mapsto \phi(u; \vartheta)$. By differentiating w.r.t.\ $\vartheta$, this is equivalent to solving
	\eqn{
	\expec[u\SS_u\e^{\vartheta u \SS_u}]=0,
	}
which in turn is equivalent to
	\eqn{
	\label{mean-SSu-tilde}
	\widetilde\expec[u\SS_u]=0,
	}
so that $\SS_u$ has mean zero under the tilted measure.
Unfortunately, \eqref{mean-SSu-tilde} turns out to be a difficult analytical problem,
and we need to resort to an approximation instead. Let us explain this in more detail now.
By the independence of the indicators $(\II_i(u))_{i\geq 2}$, we obtain
that
    \eqan{
    \lbeq{expectation-exp}
    \phi(u; \vartheta)&=\expec[\e^{\vartheta u \SS_u}]=\e^{\vartheta u(1+\betam u)}\prod_{i=2}^{\infty}\e^{-\vartheta u^2 c_i^2}
    \Big(\e^{- c_iu}
    +\e^{\vartheta c_iu}
    (1-\e^{-c_iu})\Big)\\
    &=\e^{\vartheta u(1+\betam u)} \e^{\sum_{i=2}^{\infty}f(i/u^{\tau-1};\vartheta)}\nn
    }
with (substitute $uc_i=x^{-\alpha}$)
    \eqn{
    \lbeq{def-f-tail}
    f(x;\vartheta)=\log\big(1+\e^{-x^{-\alpha}}
    (\e^{-\vartheta x^{-\alpha}}-1)\big)+\vartheta x^{-\alpha}-\vartheta x^{-2\alpha},
    }
where $\alpha=1/(\tau-1)$.
It is not hard to see that $x\mapsto f(x;\vartheta)$ is integrable at $x=0$ and at $x=\infty$ (see Lemma \ref{lem-f-integrable} below),
so we can approximate the above sum by an integral
    \eqn{
    \lbeq{sum-to-int}
    \sum_{i=2}^{\infty}f(i/u^{\tau-1};\vartheta)=u^{\tau-1}
    \int_0^{\infty}f(x;\vartheta) dx+e_\vartheta(u)
    \equiv u^{\tau-1}\Fvartheta(\vartheta)+ e_\vartheta(u),
    }
for some error term $u\mapsto e_\vartheta(u)$. For
$u$ large, the error term $u\mapsto e_\vartheta(u)$ is determined in Lemma \ref{lem-error-for-tail} below as
	\eqn{
	e_\vartheta(u)=\vartheta \left\{ u [\zeta(\alpha)-1] - u^2[\zeta(2\alpha)-1] \right\} + o_\vartheta(1),
	}
where $\zeta(\alpha),\zeta(2\alpha)$ are defined in \eqref{rz},
and where the error term converges to 0 uniformly for $\vartheta$ in compact sets bounded away from $0$.
This implies that
    \eqn{
    \phi(u; \vartheta)
    =\e^{u^{\tau-1}\Fvartheta(\vartheta)+\vartheta u(\zeta(\alpha)+(\betam-\zeta(2\alpha)+1) u)+o_\vartheta(1)}.
    }
Rather than minimizing $\phi(u; \vartheta)$ over $\vartheta$, instead we minimize the asymptotic form
appearing in its exponential $\vartheta\mapsto u^{\tau-1}\Fvartheta(\vartheta)+\vartheta u(\zeta(\alpha)+(\betam-\zeta(2\alpha)+1) u)$. Thus,
let $\thetastar_u$ be the solution of
	\eqn{
	\lbeq{thetastaru-def}
	\thetastar_u=\arg\min_{\vartheta} \Big[\Fvartheta(\vartheta)
	+\vartheta u^{2-\tau}(\zeta(\alpha)+(\betam-\zeta(2\alpha)+1) u)\Big],
	}
and let $\thetastar$ be the value of $\vartheta$ where $\vartheta\mapsto \Fvartheta(\vartheta)$ is minimal.
It is not hard to see that $I\equiv-\Fvartheta(\thetastar)>0$ and that $\thetastar$ is
\emph{unique} (see Lemma \ref{lem-unique-theta} below).
As it turns out, this choice is asymptotically equivalent to $\argmin_{\vartheta}\phi(u; \vartheta)$,
but it is analytically much more tractable. Naturally, the statement that $\thetastar$ is
asymptotically equivalent to $\argmin_{\vartheta}\phi(u; \vartheta)$ requires a proof, which can be found in
Lemma \ref{extralemma}, where we show that $\widetilde\expec[u\SS_u]=o(1)$ for $\vartheta=\thetastar_u$, and Lemma \ref{lem-approx-VP} where we show that $\thetastar_u\rightarrow \thetastar$ as $u\rightarrow \infty$.

Define $\phi(u)=\phi(u;\thetastar_u)$. The next result investigates the main term $\phi(u)$:

\begin{Proposition}[Asymptotics of main term]
\label{prop-asy-phi}
As $u\rightarrow \infty$, and with $I=-\min_{\vartheta\geq 0} \Fvartheta(\vartheta)>0$,
there exist $\kappa_{ij}\in {\mathbb R}$ such that
    \eqn{
    \lbeq{phi-asy}
    \phi(u)=\expec[\e^{\thetastar_u u \SS_u}]=\e^{-Iu^{\tau-1}+u^{\tau-1}\sum_{i+j \geq 1} \kappa_{ij}u^{-i(\tau-2)-j(\tau-3)}}(1+o(1)).
    }
\end{Proposition}
Proposition \ref{prop-asy-phi} will be proved in Section \ref{sec-main-term}.


\subsection{Properties of the process under the tilted measure}\label{proppp}
Define, for $a\in [0,1]$,
\eqn{
    \lbeq{def_I_E}
    I_{\sss E}(a)=(\tau-1)\int_0^{\infty}\Big( \frac{ \e^{\theta^* v} ( 1 - \e^{-a v} ) }{ \e^{\theta^* v}(1- \e^{-v}) + \e^{-v} } - a v\Big) \frac{ dv }{ v^{\tau - 1} }.
    }
As we see in Theorem \ref{thm-sample-path-LD}, the function $a\mapsto I_{\sss E}(a)$ will serve to describe as the
asymptotic mean of the process $a\mapsto \SS_{au}$ conditionally on $\SS_u > 0$.
It is not hard to check that
	\eqn{
    	I_{\sss E}(0)=0,\qquad \text{ and }
    	I_{\sss E}(1)=0,
    	}
the latter by definition of $\thetastar$, since $0 = \Fvartheta'(\thetastar) = I_{\sss E}(1)$ (cf.\ \eqref{der-theta-F} below). Finally,
	\eqn{
    	I_{\sss E}(a)>0 \mbox{ for every } a\in(0,1)
	}
and
	\eqn{
    	I_{\sss E}'(0)>0 \mbox{ and } I_{\sss E}'(1)<0,
	}
since $I_E$ is continously differentiable and concave on $[0,1]$ being an integral of a concave function.

From now on, we will take $\vartheta=\thetastar_u$, and we define
$\widetilde \prob=\widetilde \prob_{\thetastar_u}$ with corresponding expectation
$\widetilde \expec=\widetilde \expec_{\thetastar_u}$. In what follows, we abbreviate
$\theta=\thetastar_u$. Under this new measure, the rare event of $\SS_u$ being positive becomes quite likely,
as reflected in the following properties:

\begin{Lemma}[Expectation of $\SS_t$]
\label{lem-expec-S}
As $u\to\infty$,\\
{\rm(a)}  $\widetilde \expec[\SS_{t}] = u^{\tau-2}I_{\sss E}(t/u)+ O(1+t + t |\theta^*- \theta_u^*| u^{\tau-3})$ uniformly in $t\in [0,u]$. \\
{\rm(b)}  $\widetilde \expec[\SS_{t}-\SS_u] = u^{\tau-2}I_{\sss E}(t/u)+ O(u-t + u^{-1} + |\theta^*- \theta_u^*| u^{\tau-2})$ uniformly in $t\in [u/2,u]$.\\
{\rm(c)}  $\widetilde \expec[\SS_{t}-\SS_u] = u^{\tau-3}I_{\sss E}'(1)(t-u)(1+o(1))+ O(u^{-1})$ when $u-t=o(u)$.\\
{\rm(d)}  $u\widetilde \expec[\SS_u]=o(1)$ when $u\rightarrow \infty$.
\end{Lemma}
\medskip

The next lemma concerns the variance of the process.
Define, for $a\in [0,1]$,
\eqn{
    \lbeq{IV-def}
     I_{\sss V}(a)= (\tau-1)\int_0^{\infty} \frac{ \e^{\theta^* v}(1 - \e^{-a v} ) }{ \e^{\theta^* v}(1- \e^{-v}) + \e^{-v} } \Big(1- \frac{ \e^{\theta^* v} ( 1 - \e^{-a v} ) }{ \e^{\theta^* v}(1- \e^{-v}) + \e^{-v} } \Big) \frac{ dv }{ v^{\tau - 2} }
    }
    and
	\eqan{
	J_{\sss V}(a)&=  (\tau-1)\int_0^\infty \frac{ \e^{\theta^* v}(  \e^{-a v} - \e^{-v}) }
	{ \e^{\theta^* v}(1- \e^{-v}) + \e^{-v} } \Big(1- \frac{ \e^{\theta^* v} ( \e^{-a v} - \e^{-v}) }
	{ \e^{\theta^* v}(1- \e^{-v}) + \e^{-v} } \Big) \frac{ dv }{ v^{\tau - 2} },\\
	G_{\sss V}(a) &= (\tau-1)\int_0^\infty \frac{ \e^{2\theta^* v}( 1-\e^{-a v})(\e^{-av} - \e^{-v}) }
	{ (\e^{\theta^* v}(1- \e^{-v}) + \e^{-v})^2 }{dv \over v^{\tau-2}}. \lbeq{GV-def}
	}
Again, it is not hard to see that
	\eqn{
	\lbeq{IV-domain}
    	0<I_{\sss V}(a)<\infty \mbox{ for every } a\in(0,1], \mbox{ while } I_{\sss V}(0)=0.
	}
Similarly,
	\eqn{
	\lbeq{JV-domain}
    	0<J_{\sss V}(a)<\infty \mbox{ for every } a\in[0,1), \mbox{ while } J_{\sss V}(1)=0.
	}

\begin{Lemma}[Covariance structure of $\SS_t$]
\label{lem-var-S-rep}
As $u\to\infty$,\\
{\rm(a)}  $\widetilde{{\rm Var}}[\SS_{t}] = u^{\tau -3}I_{\sss V}(t/u) + O(1 + t |\theta^* - \theta_u^*| u^{\tau-4})$ uniformly in $t\in[0,u]$. \\
{\rm(b)}  $\widetilde{{\rm Var}}[\SS_{t}-\SS_u] = u^{\tau -3}J_{\sss V}(t/u) + O((u-t)u^{-1} + (u-t) |\theta^* - \theta_u^*| u^{\tau-4})$ uniformly in $t\in[0,u]$.\\
{\rm(c)}  $\widetilde {\rm Cov}[\SS_t,\SS_u-\SS_t] = -u^{\tau-3}G_{\sss V}(t/u)+ O((u-t)u^{-1} + (u-t) |\theta^* - \theta_u^*| u^{\tau-4})$ uniformly in $t\in[0,u]$.
\end{Lemma}	
We prove Lemma \ref{lem-expec-S} and Lemma \ref{lem-var-S-rep} in Section \ref{sec-prop-tilted}.

We complete this section by a result on the Laplace transform of the couple $(\SS_t,\SS_u)$:
\begin{Proposition}[Joint moment generating function of $(\SS_t,\SS_u)$]
\label{prop-laplace}
(a) As  $u\to \infty$,
	\eqn{
	\widetilde \expec\Big[\e^{ \lambda \frac{\SS_{t}-\widetilde \expec[\SS_t]}{\sqrt{I_{\sss V}(t/u) u^{\tau-3}}}}\Big]
	= \e^{{1\over 2}\lambda^2+ \Theta},
	}
where $|\Theta|\le o_u(1)$ as $u\rightarrow \infty$ uniformly in $t\in[u/2,u]$ and $\lambda$ in a compact set.\\
(b) Fix $\vep>0$ small.  As  $u\to \infty$, for any $\lambda_1, \lambda_2 \in \mathbb{R}$,
	\eqn{
	\widetilde \expec\Big[\e^{ \lambda_1 {\SS_{t} -\widetilde \expec[\SS_t] \over
	\sqrt{I_{\sss V}(t/u) u^{\tau-3}}} + \lambda_2  {\SS_{u} - \SS_t - \widetilde \expec[\SS_u - \SS_t]
	\over \sqrt{J_{\sss V}(t/u)u^{\tau-3}}}}\Big]
	= \e^{{1\over 2}\lambda_1^2 + {1\over 2}\lambda_2^2 -
	\lambda_1\lambda_2{G_{\sss V}(t/u)\over I_{\sss V}(t/u)J_{\sss V}(t/u)} + \Theta},
	}
where $|\Theta|\le o_u(1) + O(t^{3(3-\tau)/2})$ uniformly in $t\in[\vep,u-u^{-(\tau-5/2)}]$ and $\lambda_1, \lambda_2$ in a compact set.
\end{Proposition}
\medskip

Proposition \ref{prop-laplace} is proved in Section \ref{profff}.
By Proposition \ref{prop-laplace} and the fact that $u\widetilde \expec[\SS_u]=o(1)$
(see Lemma \ref{extralemma}), $u^{-(\tau-3)/2}\SS_u$ converges to a normal distribution with mean $0$ and variance $I_{\sss V}(1)$. We next extend this intuition by proving that the density of
$\SS_{u}$ close to zero behaves like $\left(2\pi I_{\sss V}(1)\right)^{-1/2}u^{-(\tau-3)/2}$:

\begin{Proposition}[Density of $\SS_u$ near zero]
\label{lem-dens-f}
Uniformly in  $s=o(u^{(\tau-3)/2})$, the density $\widetilde f_{\SS_u}$ of $\SS_u$ under $\widetilde \prob$ satisfies
    \eqn{
    \widetilde f_{\SS_u}(s)=B u^{-(\tau-3)/2}(1+o(1)),
    }
with $B= \left(2\pi I_{\sss V}(1)\right)^{-1/2}$ and  $I_{\sss V}(1)$ defined in  \eqref{IV-def}.
Moreover,
$\widetilde f_{\SS_t}(s)$ is uniformly bounded by a constant times $u^{-(\tau-3)/2}$ for all
$s, u$ and $t \in [u/2,u]$.
\end{Proposition}
Proposition \ref{lem-dens-f} is proved in Section \ref{profff}.

Theorems \ref{thm-tail-Su} and \ref{thm-sample-path-LD} are proved in Section \ref{sec:proofld}.
In particular, Proposition \ref{lem-dens-f} implies Theorem \ref{thm-tail-Su}, and Proposition \ref{prop-laplace} is the crucial ingredient for proving Theorem \ref{thm-sample-path-LD}.

\section{The main term: proof of Proposition \ref{prop-asy-phi}}
\label{sec-main-term}

In this section, we investigate the main term $\phi(u; \vartheta)=\expec[\e^{\vartheta u \SS_u}]$.
We want to take $\vartheta$ such that $\phi(u; \vartheta)$ is close to minimal.
Differentiating $\phi(u; \vartheta)$ with respect to $\vartheta$ suggests that we should take
$\vartheta$ such that $\expec[u \SS_u\e^{\vartheta u \SS_u}]=0$,
which is equivalent to $\widetilde\expec_\vartheta[u \SS_u]=0$. Unfortunately,
our analytical control over $\widetilde\expec_\vartheta[u \SS_u]$ is too limited to make this choice work, so instead
we optimize the asymptotic expression \eqref{expectation-exp} for $\phi(u; \vartheta)$ instead.
To this end, the main result in this section is Lemma \ref{lem-error-for-tail},
which sets the stage for the proof of Proposition \ref{prop-asy-phi}.

We start by proving properties of the function $f$ defined
in \refeq{def-f-tail}.

\begin{Lemma}[Integrability of $f(\cdot;\vartheta)$]
\label{lem-f-integrable}
Fix $\vartheta>-1$. The function $x\mapsto f(x;\vartheta)$ with $f$ as in \eqref{def-f-tail} is integrable at $x=0$ and at $x=\infty$.
\end{Lemma}

\proof For $x \downarrow 0^+$ and $\vartheta>-1$, the first term of $f(x;\vartheta)$ approaches zero and the second and third term are integrable at $x=0^+$. The case where $x \rightarrow \infty$ requires to consider the conjunction of all three terms. Let $y = x^{-\alpha}$, so that we have to consider integrability at $y=0^+$. We can use Taylor approximation to obtain
    \eqan{
    \lbeq{integrability}
    f(x;\vartheta)
    &= \log\big(1+\e^{-y} (\e^{-\vartheta y}-1)\big)+\vartheta y-\vartheta y^2
    = y^3 \frac{1}{2} \vartheta (\vartheta-1) + O(y^4)\\
    &= x^{-3\alpha} \frac{1}{2} \vartheta (\vartheta-1) + O\!\left( x^{-4\alpha} \right),\nn
    }
which is integrable for $x \rightarrow \infty$ since $-3\alpha \in (-3/2,-1)$.
\qed
\medskip

We continue with a general result allowing us to replace sums
by integrals with a good control over the error term.

\begin{Lemma}[Approximating sums by integrals]
\label{lem-taylor}
Let $g\colon\mathbb{R}_+\to \mathbb{R}$ be a differentiable function such that there exist $\gamma>-1$, and $a,b\ge 0$ satisfying $|g'(y)|\le ay^{\gamma}\e^{-by}$ for all $y> 0$. Then, for any $\alpha >0$, there exist $c=c(\alpha,\gamma)$ {\rm(}which does not depend on $g${\rm)} such that
	$$
	\Big|\sum_{i=2}^{\infty} \big[g(u{i^{-\alpha}}) - \int_{i}^{i+1} g(ux^{-\alpha})dx\big]\Big| \le
	c a  \min(u,b^{-1})^{1+\gamma}.
	$$
\end{Lemma}
\medskip


\proof  Let $h(x)=g(ux^{-\alpha})$. By the Taylor approximation $|h(x)-h(i)|\le \sup_{y\in [i,i+1]} |h'(y)| (x-i)$ for any $x\in[i,i+1]$, for any $i\ge 2$,
	$$
	\Big| g(u{i^{-\alpha}}) - \int_{i}^{i+1} g(ux^{-\alpha})dx \Big|
	\le {1\over 2}\sup_{x\in[i,i+1]} |h'(x)|.
	$$
By assumption, we know that $\sup_{x\in[i,i+1]} |h'(x)|$ is less than $a \alpha u^{1+\gamma}i^{-\alpha (1+\gamma)-1}\e^{-bu(i+1)^{-\alpha}}$. This yields that
	\eqn{\label{eq:taylor}
	\sum_{i\ge 2} \Big|g(u{i^{-\alpha}}) - \int_{i}^{i+1} g(ux^{-\alpha})dx\Big|
	\le
	a \alpha u^{1+\gamma}\sum_{i\ge 2}i^{-\alpha (1+\gamma)-1}\e^{-bu(i+1)^{-\alpha}}.
	}
For any $x\in[i,i+1]$ and $i\ge 2$,
	$$
	i^{-\alpha (1+\gamma)-1}\e^{-bu(i+1)^{-\alpha}} \le 2^{\alpha (1+\gamma)+1}
	x^{-\alpha (1+\gamma)-1}\e^{-bu x^{-\alpha}}.
	$$
Hence,
	$$
	\sum_{i\ge 2 }i^{-\alpha (1+\gamma)-1}\e^{-bu(i+1)^{-\alpha}}
	\le
	c(\alpha,\gamma) \int_{2}^{\infty} x^{-\alpha (1+\gamma)-1}\e^{-bu x^{-\alpha}}dx,
	$$
which is less than $c(\alpha,\gamma) (bu)^{-(1+\gamma)}$. Since $\gamma+1>0$, we have as well $\sum_{i\ge 2}i^{-\alpha (1+\gamma)-1}\e^{-bu(i+1)^{-\alpha}} \le c(\alpha,\gamma)$. It follows that
	$$
	\sum_{i\ge 2 }i^{-\alpha (1+\gamma)-1}\e^{-bu(i+1)^{-\alpha}}
	\le
	c(\alpha,\gamma) \min(1,(bu)^{-(1+\gamma)}).
	$$
Then (\ref{eq:taylor}) completes the proof.
\qed

\begin{Corollary}[Replacing sums by integrals in general]
\label{cor-exp-mom-cs}
For every $a\in \R, a > \tau-1$ and $b>0$, there exists a constant $c(a,b)$ such that
	\eqn{
	\sum_{i=2}^{\infty} c_i^a \e^{-b c_i u}=c(a,b) u^{\tau-a-1}(1+o(1)).
	}
\end{Corollary}

\proof This follows directly from Lemma \ref{lem-taylor} with $g(y)=y^a \e^{-by}.$
\qed
\medskip

We next investigate the error in replacing the sum over $i$ of $f(i/u^{\tau-1};\vartheta)$
by the integral in \eqref{sum-to-int}, using similar ideas as in Lemma \ref{lem-taylor}
above.

\begin{Lemma}[Error in replacing sum by integral]
\label{lem-error-for-tail}
The error term $e_\vartheta(u)$ in \eqref{sum-to-int} satisfies
    \eqan{
    \lbeq{F1-and-F2}
    e_\vartheta(u)
    &= \vartheta \left\{  u \left[ \zeta(\alpha) -  1 \right] - u^2 \left[ \zeta(2\alpha) - 1  \right] \right\} + o_\vartheta(1)
    }
where $o_\vartheta(1)$ depends on $\vartheta$ and satisfies that, uniformly in $\vartheta\in (-1+\vep, 1/\vep)$,
    \eqn{
    \lbeq{order-vartheta}
    |o_\vartheta(1) | \leq C_\vep u^{-(\tau-1)}
    }
for $u$ large enough. Further,
	\eqn{
	\lbeq{deriv-e-bd}
	e_{\vartheta}'(u)
	=\sum_{i\geq 2} \frac{\partial}{\partial \vartheta} f(i/u^{\tau-1}; \vartheta)
	-u^{\tau-1}\int_{0}^{\infty} \frac{\partial}{\partial \vartheta} f(x; \vartheta)dx
	=u (\zeta(\alpha)-1) - u^2 (\zeta(2\alpha)-1)+ o_\vartheta(1),
	}
where $o_\vartheta(1)$ satisfies \eqref{order-vartheta} as well.
\end{Lemma}
\medskip

Note that \refeq{SS-mean-ass} follows by taking $\vartheta=0$ in \refeq{deriv-e-bd} and using \eqref{expectation-exp} and \eqref{sum-to-int} to get
	\[
	\expec[\SS_u]=\frac{1}{u} \left.
	\frac{\partial}{\partial \vartheta} \right|_{\vartheta=0}
	\log\big( \expec[\e^{\vartheta u \SS_u}] \big) = u^{\tau-2}\Fvartheta'(0)+(1+\betam u)+e_{0}'(u)/u.
	\]
The precise form of $\Fvartheta'(0)$ can be found in \eqref{Fvartheta-der-0} below.

%
%
%

\proof For $\vartheta>-1+\vep$, we split
    \eqn{
    f(x;\vartheta)
    =\log\big(1+\e^{-x^{-\alpha}}
    (\e^{-\vartheta x^{-\alpha}}-1)\big)+\vartheta x^{-\alpha}-\vartheta x^{-2\alpha}
    \equiv f_1(x;\vartheta)+f_2(x;\vartheta)+f_3(x;\vartheta)
    }
(cf.~\eqref{def-f-tail}). Remember that $\tau \in (3,4)$ and thus $\alpha \in (1/3,1/2)$.

In what follows, fix $u$ arbitrarily large. By Lemma \ref{lem-f-integrable}, for every $\varepsilon>0$ small, we can choose $M=M(u)>0$ large such that
    \eqn{
    \sum_{i \geq Mu^{\tau-1}} \Big| f\!\Big( \frac{i}{u^{\tau-1}} ; \vartheta \Big) \Big| + u^{\tau-1} \int_M^\infty \left| f(x;\vartheta) \right| dx < \varepsilon.
    }
It remains to estimate the difference
    \eqn{
     \sum_{2 \leq i \leq Mu^{\tau-1}-1} f\!\Big( \frac{i}{u^{\tau-1}} ; \vartheta \Big) - u^{\tau-1} \int_0^M f(x;\vartheta) dx
    }
for $M$ arbitrarily large.

We analyse the second and third term in the definition of $f(x;\vartheta)$ for $i \leq M u^{\tau-1}-1$ respectively $x \leq M$ first. Observe that for arbitrary $\beta \in (0,1)$ (later to be set equal to $\alpha$, respectively $2\alpha$)
    \eqan{
    & \sum_{2 \leq i \leq Mu^{\tau-1}-1} \left( \frac{i}{u^{\tau-1}} \right)^{-\beta} - u^{\tau-1} \int_0^M x^{-\beta} dx \\
    &= u^{(\tau-1) \beta} \sum_{2 \leq i \leq Mu^{\tau-1}-1} i^{-\beta} - u^{(\tau-1)\beta} \int_1^{M u^{\tau-1}} x^{-\beta} dx - u^{\tau-1} \int_0^{\frac{1}{u^{\tau-1}}} x^{-\beta} dx \nn\\
    &= u^{(\tau-1)\beta} \sum_{2 \leq i \leq Mu^{\tau-1}-1} i^{-\beta} - u^{(\tau-1)\beta} \int_1^{M u^{\tau-1}} x^{-\beta} dx - u^{\tau-1} \frac{1}{1-\beta} u^{-(\tau-1)(1-\beta)} \nn\\
    &= u^{(\tau-1)\beta} \Big( \sum_{1 \leq i \leq Mu^{\tau-1}-1} i^{-\beta} - \int_1^{M u^{\tau-1}} x^{-\beta} dx \Big) - u^{(\tau-1)\beta} \Big( \frac{1}{1-\beta} + 1 \Big). \nn
    }
Here,
    \eqn{
    \sum_{1 \leq i \leq Mu^{\tau-1}-1} i^{-\beta} - \int_1^{M u^{\tau-1}} x^{-\beta} dx = \sum_{1 \leq i \leq Mu^{\tau-1}-1} \Big( i^{-\beta} - \int_i^{i+1} x^{-\beta} dx \Big) \equiv C(M,\beta)
    }
with
    \eqn{
    C(M,\beta) \uparrow C(\beta) \equiv \sum_{i\geq 1} \left( i^{-\beta} - \int_i^{i+1} x^{-\beta} dx \right), \quad M \rightarrow \infty
    }
satisfying $0<C(\beta)<\infty$. Indeed, $x^{-\beta}$ is a strictly decreasing function and therefore
    \eqn{
    0 < i^{-\beta} - \int_i^{i+1} x^{-\beta} dx < i^{-\beta} - (i+1)^{-\beta} \leq (-\beta) i^{-\beta-1} (-1),
    }
and hence
    \eqn{
    \lbeq{upper-bound-on-C-beta}
    0 < \sum_{i \geq 1} \left[ i^{-\beta} - \int_i^{i+1} x^{-\beta} dx \right] \leq \beta \sum_{i \geq 1} i^{-\beta-1} < \infty
    }
for all $\beta \in (0,1)$. For $\tau \in (3,4)$, we have $\alpha \in (1/3,1/2)$ and $2\alpha \in (2/3,1)$, so we can apply the above result to both the second term ($\beta=\alpha$) and the third term ($\beta=2\alpha$) to obtain
    \eqan{
    & \lim_{M \rightarrow \infty} \left\{ \sum_{2 \leq i \leq Mu^{\tau-1}-1} (f_2+f_3)\!\left( \frac{i}{u^{\tau-1}} ; \vartheta \right) - u^{\tau-1} \int_0^M (f_2+f_3)(x;\vartheta) dx \right\} \\
    &= \vartheta \left\{  u \left[ C(\alpha) - \left( \frac{1}{1-\alpha} + 1 \right) \right] - u^2 \left[ C(2\alpha) - \left( \frac{1}{1-2\alpha} + 1 \right) \right] \right\}, \nn \\
    &= \vartheta \left\{  u \left[ \zeta(\alpha) -  1 \right] - u^2 \left[ \zeta(2\alpha) -  1 \right] \right\}, \nn
    }
where we have used $(\tau-1)\alpha=1$, and where the last identity follows from \eqref{rz}.

It remains to analyse the contribution due to $f_1(x;\vartheta)=\log\big(1+\e^{-x^{-\alpha}}
    (\e^{-\vartheta x^{-\alpha}}-1)\big)$. Observe also that,
$\lim_{x \downarrow 0} f_1(x;\vartheta) = 0$ since $\vartheta>-1$.
We first calculate the first two derivatives of $f_1(x;\vartheta)$ with respect to $x$.
We note that $f_1(x;\vartheta)=g(x^{-1/\alpha}),$ where $g(x)=\log\big(1+\e^{-x}
 (\e^{-\vartheta x}-1)\big)$. Therefore,
	\eqan{
	\frac{\partial}{\partial x}f_1(x;\vartheta)&=-(1/\alpha) x^{-(1/\alpha+1)} g'(x^{-1/\alpha}),\\
	\frac{\partial^2}{\partial x^2}f_1(x;\vartheta)
	&=(1/\alpha)(1/\alpha+1)x^{-(1/\alpha+2)} g'(x^{-1/\alpha})
	+(1/\alpha+2) x^{-(2/\alpha+2)} g''(x^{-1/\alpha}),
	}
with
	\eqn{
	g'(x)=-\frac{\e^{-x}
	\left( (\vartheta+1) \e^{-\vartheta x} - 1 \right)}{1+\e^{-x}
    	(\e^{-\vartheta x}-1)},
	\qquad
	g''(x)=[g'(x)]^2+\frac{\e^{-x}
	\left( (\vartheta+1)^2 \e^{-\vartheta x} - 1 \right)}{1+\e^{-x}
    	(\e^{-\vartheta x}-1)}
	}
In particular, $|g''(x)|\leq C\e^{-x}$ uniformly in $x\geq 0$, so that
$|\frac{\partial^2}{\partial x^2}f_1(x;\vartheta)|$ is integrable, and is uniformly integrable for
$\vartheta>-1+\vep$ for any $\vep>0$.

	
Rewrite
    \eqan{
    \lbeq{second-diff}
    & \sum_{2 \leq i \leq Mu^{\tau-1}-1} f_1(i/u^{\tau-1}; \vartheta) - u^{\tau-1} \int_0^M f_1(x;\vartheta) dx \\
    & =
	\sum_{1 \leq i \leq Mu^{\tau-1}-1} \left[ f_1(i/u^{\tau-1}; \vartheta) -  \int_{i-1/2}^{i+1/2}
	f_1(x/u^{\tau-1}; \vartheta) dx \right]\nn\\
	&\qquad
	- f_1(1/u^{\tau-1}; \vartheta)+\int_0^{1/2} f_1(x/u^{\tau-1}; \vartheta) dx
	+ \int_{M-1/2}^{M} f_1(x/u^{\tau-1}; \vartheta) dx.\nn
    }
We first identify the error terms. Using $|\log\big(1-h\big)|\leq |h|/(1-h)$ for $|h|>1$
we obtain
    	\eqn{
    	|f_1(1/u^{\tau-1}; \vartheta)|=
	\big| \log\big(1+\e^{-u} (\e^{-\vartheta u}-1)\big) \big| \leq \frac{\e^{-[(1+\vartheta)\wedge 1)]u}}{1-\e^{-[(1+\vartheta)\wedge 1]u}}.
    	}
The term $\int_0^{1/2} f_1(x/u^{\tau-1}; \vartheta) dx$ can be seen to obey a similar bound.
Finally, since $\lim_{y\rightarrow \infty} f_1(y; \vartheta)=0$, the term $\int_{M-1/2}^{M} f_1(x/u^{\tau-1}; \vartheta) dx$ can be made arbitrarily small by taking $M$ large.

By quadratic Taylor approximation, for $x\in [i-1/2, i+1/2]$,
	\eqn{
	f_1(x/u^{\tau-1}; \vartheta) =f_1(i/u^{\tau-1}; \vartheta) +u^{-(\tau-1)}(x-i)\frac{\partial}{\partial x}f_1(i/u^{\tau-1}; \vartheta) +\frac{u^{-2(\tau-1)}}{2}(x-i)^2\frac{\partial^2}{\partial x^2}f_1(\xi/u^{\tau-1}; \vartheta),
	}
for some $\xi\in[i-1/2,i+1/2]$. Note that the integral over $x\in [i-1/2,i+1/2]$
of the first term equals $f_1(i/u^{\tau-1}; \vartheta)$ and of the second term equals 0.
Therefore,
    \eqan{
    &\Big| \sum_{0 \leq i \leq Mu^{\tau-1}-1} \Big[ f_1(i/u^{\tau-1}; \vartheta) -  \int_{i-1/2}^{i+1/2} f_1(x/u^{\tau-1}; \vartheta) dx\Big] \Big|\\
     &\leq \frac{1}{12}\sum_{i \geq 0} \sup_{\xi \in [i-1/2,i+1/2]}\Big|\frac{\partial^2}{\partial x^2}f_1(x/u^{\tau-1};\vartheta)\Big| u^{-2(\tau-1)}\le c(\vartheta) u^{-(\tau-1)}, \nn
    }
where we use that, as $u\to\infty$, by a Riemann sum approximation and the fact that $|\frac{\partial^2}{\partial x^2}f_1(x/u^{\tau-1};\vartheta)|$
is integrable,
	\eqn{
	u^{-(\tau-1)}\sum_{i \geq 1} \sup_{\xi \in [i-1/2,i+1/2]}
	|\frac{\partial^2}{\partial x^2}f_1(x/u^{\tau-1};\vartheta)|
	\rightarrow \int_0^{\infty}|\frac{\partial^2}{\partial x^2}f_1(x;\vartheta)|dx.
	}
The claim now follows after collecting terms and taking $M \rightarrow \infty$.

The proof for $\frac{\partial}{\partial \vartheta}f(x;\vartheta)$ is identical, now using that
	\eqn{
    	\frac{\partial}{\partial \vartheta}f(x;\vartheta)
	=x^{-\alpha}\left[1-x^{-\alpha}-\frac{1}{1+\e^{\vartheta x^{-\alpha} }(\e^{x^{-\alpha}}-1)}\right]
	=\frac{\partial}{\partial \vartheta}f_1(x;\vartheta)+f_2(x;\vartheta)/\vartheta+f_3(x;\vartheta)/\vartheta.
    	}
The sums of $f_2(i/u^{\tau-1};\vartheta)/\vartheta+f_3(i/u^{\tau-1};\vartheta)/\vartheta$ give rise to $u (\zeta(\alpha)-1) - u^2 (\zeta(2\alpha)-1)$ as in the above argument, while the sum of $\frac{\partial}{\partial \vartheta}f_1(i/u^{\tau-1};\vartheta)$ is $u^{\tau-1} \int_0^{\infty} \frac{\partial}{\partial \vartheta}f_1(x;\vartheta)dx+o_\vartheta(1)$
as for $f_1(x;\vartheta)$.
\qed
\bigskip

We next investigate the properties of $\Fvartheta(\vartheta)=\int_0^\infty f(x;\vartheta)dx$:

\begin{Lemma}[Uniqueness of limiting variational problem]
\label{lem-unique-theta}
Let $\thetastar$ be the value of $\vartheta$ where $\vartheta\mapsto \Fvartheta(\vartheta)$ is minimal. Then $\thetastar$ is the unique positive solution to $\Fvartheta'(\vartheta)=0$ and $I=-\Fvartheta(\thetastar)>0$.
\end{Lemma}

\proof We have that
    \eqn{
    \lbeq{der-theta-f}
    \frac{\partial}{\partial \vartheta}f(x;\vartheta)=x^{-\alpha}\left[1-x^{-\alpha}-\frac{1}{1+\e^{\vartheta x^{-\alpha} }(\e^{x^{-\alpha}}-1)}\right]
    }
and
    \eqn{
    \frac{\partial^2}{\partial^2 \vartheta}f(x;\vartheta)=\frac{\e^{\vartheta x^{-\alpha}}(\e^{x^{-\alpha}}-1)x^{-2\alpha}}{\left[1+\e^{\vartheta x^{-\alpha} }(\e^{x^{-\alpha}}-1)\right]^2}.
    }
Differentiation under the integral sign, using dominated convergence, yields
    \eqn{
    \lbeq{der-theta-F}
    \Fvartheta'(\vartheta)=\int_0^\infty \frac{\partial}{\partial \vartheta}f(x;\vartheta)dx,
	\qquad\quad \Fvartheta''(\vartheta)=\int_0^\infty \frac{\partial^2}{\partial^2 \vartheta}f(x;\vartheta)dx.
    }
%
%
Observe that $ \frac{\partial^2}{\partial^2 \vartheta}f(x;\vartheta)>0$ for all $\vartheta\in\mathbb{R}$, as long as $x>0$, and hence that $\Fvartheta''(\vartheta)>0$ for all $\vartheta\in\mathbb{R}$. Now, since $\Fvartheta(0)=0$, $\Fvartheta'(\vartheta)\rightarrow \infty$ as $\vartheta\rightarrow \infty$ (as we will show below) and
    \eqn{
    \lbeq{Fvartheta-der-0}
    \Fvartheta'(0)=\int_0^\infty x^{-\alpha}[1-x^{-\alpha}-\e^{-x^{-\alpha}}]dx<0,
    }
the assertion follows.

To prove that $\Fvartheta'(\vartheta)\rightarrow \infty$ as $\vartheta\rightarrow \infty$, we bound $\Fvartheta'(\vartheta)$ from below. Using \eqref{der-theta-F} and \eqref{der-theta-f} together with the substitution $y=x^{-\alpha}$,
    \eqn{
    \Fvartheta'(\vartheta) = \frac{1}{\alpha} \int_0^\infty \left[ 1-y-\frac{1}{1+\e^{\vartheta y} ( \e^y-1 ) } \right] y^{-1/\alpha} dy.
    }
We split the integral, depending on whether $y\in [0,1/2]$ or not.
For $y>1/2$, the absolute value of the integrand is uniformly bounded by $cy^{1-1/\alpha}$,
which is integrable on $[1/2,\infty)$ since $1-1/\alpha=2-\tau<-1$. For $y\in[0,1/2]$, and for $\vartheta>0$ sufficiently large, the integrand is
non-negative, increasing in $\vartheta$ and converges pointwise to $(1-y)y^{-1/\alpha}$. Therefore, by monotone convergence
	\eqn{
	\lim_{\vartheta\rightarrow \infty}\int_0^{1/2}
	\left[ 1-y-\frac{1}{1+\e^{\vartheta y} ( \e^y-1 ) } \right] y^{-1/\alpha} dy
	=\int_0^{1/2}[1-y] y^{-1/\alpha} dy=\infty,
	}
since $1/\alpha=\tau-1>1$. We conclude that
	\eqn{
	\label{F'-infty}
	\lim_{\vartheta\rightarrow \infty} \Fvartheta'(\vartheta)=\infty.
	}
The claim now follows.
\qed
\medskip

We continue to investigate the approximate variational problem as formulated
in \refeq{thetastaru-def}, and prove Proposition \ref{prop-asy-phi}.

\begin{Lemma}[Expansion for the maximizer $\thetastar_u$]
\label{lem-approx-VP}
Define
	\eqn{
	\lbeq{def-vep-u}
	\vep_u=u^{2-\tau}(\zeta(\alpha)+(\betam-\zeta(2\alpha)+1) u).
	}
Then, there exist $\kappa_i'\in {\mathbb R}$ such that for each $m\geq 0$
	\eqn{
	\lbeq{thetastar-asymp}
	\thetastar_u=\thetastar+\sum_{i=1}^m \kappa_i' \vep_u^i
	+O(\vep_u^{m+1}).
	}
\end{Lemma}


\proof The function $\vartheta\mapsto \Fvartheta(\vartheta)$ is differentiable.
Therefore, the minimizer $\thetastar_u$ of $\vartheta\mapsto \Fvartheta(\vartheta)+\vartheta\vep_u$
(cf. \eqref{thetastaru-def}) satisfies
	\eqn{
	\Fvartheta'(\thetastar_u)=-\vep_u.
	}
Clearly, $\vep_u\rightarrow 0$ as $u\rightarrow \infty$, and the above is an implicit
equation for $\thetastar_u$. We define $\Fvartheta^{\sss(r)}(\vartheta)$ to be the $r$-fold derivative of
$\Fvartheta$ with respect to $\vartheta$, and we let $f^{\sss(r)}(x;\vartheta)$ be the
$r$-fold derivative of $\vartheta\mapsto f(x;\vartheta)$ with respect to $\vartheta$,
where we recall that $x\mapsto f(x;\vartheta)=\log\big(1+\e^{-x^{-\alpha}}
(\e^{-\vartheta x^{-\alpha}}-1)\big)+\vartheta x^{-\alpha}-\vartheta x^{-2\alpha}$ as defined
in \refeq{def-f-tail}. Then,
	\eqn{
	\Fvartheta^{\sss(r)}(\vartheta)=\int_0^{\infty} f^{\sss(r)}(x;\vartheta)dx.
	}
We compute that
	\eqn{
	\lbeq{Init-IH-fr-1}
	f^{\sss(1)}(x;\vartheta)=\frac{-x^{-\alpha}\e^{-x^{-\alpha}}
	\e^{-\vartheta x^{-\alpha}}}{1+\e^{-x^{-\alpha}}
	(\e^{-\vartheta x^{-\alpha}}-1)}+x^{-\alpha}-x^{-2\alpha}
	=x^{-\alpha}\frac{1-\e^{-x^{-\alpha}}}{1+\e^{-x^{-\alpha}}
	(\e^{-\vartheta x^{-\alpha}}-1)}-x^{-2\alpha},
	}
and
	\eqan{
	\lbeq{Init-IH-fr}
	f^{\sss(2)}(x;\vartheta)&=x^{-2\alpha}[1-\e^{-x^{-\alpha}}]
	\frac{\e^{-x^{-\alpha}}\e^{-\vartheta x^{-\alpha}}}{[1+\e^{-x^{-\alpha}}
	(\e^{-\vartheta x^{-\alpha}}-1)]^2}.
	}
so that, in particular, $f^{\sss(2)}(x;\vartheta)>0$ for every $x, \vartheta>0$.
The latter explains why $\Fvartheta^{\sss(2)}(\vartheta)>0$ for every $\vartheta>0$.

{
We start by checking that $\vartheta \mapsto \Fvartheta(\vartheta)$ is infinitely often differentiable. {Recall \eqref{Init-IH-fr-1} and rewrite \eqref{Init-IH-fr} to
	\eqn{
	\lbeq{Init-IH-fr-2}
	f^{\sss(2)}(x;\vartheta)
	=x^{-2\alpha}\frac{-[1-\e^{-x^{-\alpha}}]^2}{[1+\e^{-x^{-\alpha}}
	(\e^{-\vartheta x^{-\alpha}}-1)]^2}
	+x^{-2\alpha}\frac{1-\e^{-x^{-\alpha}}}{1+\e^{-x^{-\alpha}}
	(\e^{-\vartheta x^{-\alpha}}-1)}.
	}}
We prove, by induction, that there exist integers $a_{r,i}$, for $i=1, \ldots, r$, such that,
for all $r\geq 2$,
	\eqn{
	\lbeq{IH-fr}
	f^{\sss(r)}(x;\vartheta)
	=x^{-r\alpha}\sum_{i=1}^r a_{r,i}\frac{[1-\e^{-x^{-\alpha}}]^i}{[1+\e^{-x^{-\alpha}}
	(\e^{-\vartheta x^{-\alpha}}-1)]^i}.
	}
We use \refeq{Init-IH-fr-2} to initialize the induction hypothesis in \refeq{IH-fr}
for $r=2$, with $a_{2,1}=-1$ and $a_{2,2}=1$. We compute that the derivative of
$\vartheta\mapsto [1-\e^{-x^{-\alpha}}]^i/[1+\e^{-x^{-\alpha}}(\e^{-\vartheta x^{-\alpha}}-1)]^i$ equals
	\eqn{
	\lbeq{Istep}
	\frac{i x^{-\alpha}[1-\e^{-x^{-\alpha}}]^i\e^{-x^{-\alpha}}
	\e^{-\vartheta x^{-\alpha}}}{[1+\e^{-x^{-\alpha}}
	(\e^{-\vartheta x^{-\alpha}}-1)]^{i+1}}
	=\frac{-i x^{-\alpha}[1-\e^{-x^{-\alpha}}]^{i+1}}{[1+\e^{-x^{-\alpha}}
	(\e^{-\vartheta x^{-\alpha}}-1)]^{i+1}}+\frac{i x^{-\alpha}[1-\e^{-x^{-\alpha}}]^{i}}{[1+\e^{-x^{-\alpha}}
	(\e^{-\vartheta x^{-\alpha}}-1)]^{i}}.
	}
We now check integrability for $x \downarrow 0^+$ respectively $x \rightarrow \infty$ for $r \geq 1$ arbitrary.
\medskip

\noindent
{\it (1) Case $r \in \{1,2\}$.}
As the denominator in \eqref{Init-IH-fr-1} respectively \eqref{Init-IH-fr-2} is uniformly greater $0$ and as $\alpha \in (1/3,1/2)$, integrability at $x \downarrow 0^+$ follows from the integrability of $x^{-\alpha}$ and $x^{-2\alpha}$. For $x \rightarrow \infty$, use that $1-\e^{-x^{-\alpha}} \sim x^{-\alpha}$ and reason as in \eqref{integrability} to obtain the claim.
\medskip

\noindent
{\it (2) Case $r \geq 3$.}
For $x \downarrow 0^+$, use \eqref{IH-fr} and \eqref{Istep} to see that
	\eqan{
	|f^{\sss(r)}(x;\vartheta)|
	&\leq \sum_{i=1}^{r-1} |a_{r-1,i}| \Big|\frac{d}{d\vartheta} \frac{[1-\e^{-x^{-\alpha}}]^i}{[1+\e^{-x^{-\alpha}}(\e^{-\vartheta x^{-\alpha}}-1)]^i}\Big|\\
	&
	\leq x^{-(r-1)\alpha}\sum_{i=1}^{r-1} |a_{r-1,i}| \frac{i x^{-\alpha}[1-\e^{-x^{-\alpha}}]^i\e^{-x^{-\alpha}}
	\e^{-\vartheta x^{-\alpha}}}{[1+\e^{-x^{-\alpha}}
	(\e^{-\vartheta x^{-\alpha}}-1)]^{i+1}}
	\leq c x^{-(r-1)\alpha} \e^{-x^{-\alpha}},\nn
	}
which is integrable for $x \downarrow 0^+$. For $x \rightarrow \infty$, use $| f^{\sss(r)}(x;\vartheta) | \leq c x^{-r\alpha}$ and $r\alpha \geq 3\alpha > 1$ in \eqref{IH-fr} to conclude integrability.
}

By \eqref{F'-infty}, $\Fvartheta^{\sss(1)}(\vartheta)\rightarrow \infty$  when $\vartheta\rightarrow \infty$.
Since $\Fvartheta^{\sss(1)}(0)<0$, $\Fvartheta^{\sss(1)}(\vartheta)\rightarrow \infty$
when $\vartheta\rightarrow \infty$ and $\Fvartheta^{\sss(2)}(\vartheta)>0$, the equation
$\Fvartheta^{\sss(1)}(\vartheta)=-\vep_u$ has a unique solution.

Let $\theta(\vep)$ be the solution $\vartheta$ to $\Fvartheta^{\sss(1)}(\vartheta)=-\vep$.
Since $\Fvartheta^{\sss(2)}(\vartheta)>0$ for every $\vartheta>0$ and we have shown that
$\vartheta \mapsto \Fvartheta^{\sss(1)}(\vartheta)$ is infinitely often differentiable, the implicit function theorem implies that $\vep\mapsto \theta(\vep)$ is
infinitely often differentiable as well in a neighborhood of $0$.
As a result, a Taylor expansion of $\vep\mapsto \theta(\vep)$ around $\vep=0$ yields that
for each $m\geq 0$ there exist $\kappa_i'\in {\mathbb R}$ such that
	\eqn{
	\theta(\vep)=\thetastar+\sum_{i=1}^m \kappa_i' \vep^i
	+O(\vep^{m+1}).
	}
Applying this identity to $\vep=-\vep_u$, while observing that $\vep_u \rightarrow 0$ for $u \rightarrow \infty$, we arrive at \refeq{thetastar-asymp}.
\qed
\medskip

Now we are ready to complete the proof of Proposition \ref{prop-asy-phi}:\\
{\it Proof of Proposition \ref{prop-asy-phi}.} By construction and Lemma \ref{lem-error-for-tail},
	\eqn{
	\phi(u)=\phi(u; \thetastar_u)
	=\e^{u^{\tau-1}[\Fvartheta(\thetastar_u)+\thetastar_u \vep_u] + o(1)}.
	}
Since $\vartheta\mapsto \Fvartheta(\vartheta)$ is infinitely differentiable,
a Taylor expansion yields
	\eqn{
	\Fvartheta(\thetastar_u)=\Fvartheta(\thetastar)+\sum_{r=1}^{m} \frac{\Fvartheta^{\sss(r)}(\thetastar)}{r!}(\thetastar_u-\thetastar)^r
	+O(|\thetastar_u-\thetastar|^{m+1}).
	}
By Lemma \ref{lem-approx-VP}, $|\thetastar-\thetastar_u|^{m+1}=O(|\vep_u|^{m+1})=o(u^{-(\tau-1)})$ if $m$ satisfies $(3-\tau)(m+1) < -(\tau-1)$, so that $u^{\tau-1}|\thetastar-\thetastar_u|^{m+1}=o(1)$. By \eqref{thetastar-asymp},
$\thetastar_u=\thetastar+\sum_{i=1}^m \kappa_i' \vep_u^i+O(\vep_u^{m+1})$.
Note that using \eqref{def-vep-u}
	\eqn{
	\vep_u^q=\sum_{i+j=q} {q \choose i} \zeta(\alpha)^i (\betam-\zeta(2\alpha)+1)^j u^{i(2-\tau)+j(3-\tau)}.
	}
Rearranging sums, we obtain that there exist $\kappa_{ij}$ such that
	\eqn{
	u^{\tau-1}[\Fvartheta(\thetastar_u)+\thetastar_u \vep_u]
	= u^{\tau-1}\Fvartheta(\thetastar)+u^{\tau-1}\sum_{i,j} \kappa_{ij} u^{i(2-\tau)+j(3-\tau)}+o(1).
	}
Since $I=-\Fvartheta(\thetastar)$, this completes the proof of Proposition \ref{prop-asy-phi}.
\qed


\section{Properties of process under tilted measure}
\label{sec-prop-tilted}

Fix $u\ge 0$. In \eqref{tildeP-def}, we defined the measure $\tilde \prob$ with Radon-Nikodym derivative $\e^{\vartheta u \SS_u} /\phi(u)$ with respect to $\prob$, where
$\vartheta=\theta_u^*$ as in \eqref{thetastaru-def}.
In particular, we stress that $\widetilde \prob$ depends on $u$. This section is devoted to the study of $(\SS_t,\,t\in [0,u])$ under $\widetilde \prob$. We derive asymptotics of
$\widetilde\expec[\SS_{t}]$, the variance of $\SS_{t}$,  and the covariance of $\SS_{t}$ and $\SS_{u}-\SS_t$,
 for all $t\in [0,u]$ that allow us to prove Lemma \ref{lem-expec-S},
and Lemma \ref{lem-var-S-rep}.

As before, and throughout the remainder of this paper, we fix $\vartheta=\theta_u^*$.
We start by proving that $u\widetilde\expec[\SS_{u}]$ vanishes as $u\rightarrow \infty$.
Intuitively, this means that we have chosen $\vartheta=\theta_u^*$ asymptotically correct:

\begin{Lemma}[Mean under tilted measure]
\label{extralemma}
As $u\rightarrow \infty$,
	\eqan{
	u\widetilde\expec[\SS_{u}]=o(1).
	}
\end{Lemma}

\proof
Note that, by \eqref{phi-def}, \eqref{expectation-exp} and \eqref{sum-to-int}, together with \refeq{deriv-e-bd} in Lemma \ref{lem-error-for-tail}
	\eqan{
	u\widetilde \expec[\SS_{u}]
	&= \frac{1}{\phi(u;\theta_u^*)} \left. \frac{\partial}{\partial \vartheta} \right|_{\vartheta=\theta_u^*} \expec[ \e^{\vartheta u \SS_u} ]
	=\Big[ u(1+\betam u) + \sum_{i\geq 2} \left. \frac{\partial}{\partial \vartheta} \right|_{\vartheta=\theta_u^*} f(i/u^{\tau-1}; \vartheta) \Big]\\
	&= \Big[ u(1+\betam u) + u^{\tau-1}\int_{0}^{\infty} \left. \frac{\partial}{\partial \vartheta} \right|_{\vartheta=\theta_u^*} f(x; \vartheta)dx +e_{\theta_u^*}'(u) \Big] \nn\\
	&=\Big[ u^{\tau-1}\Fvartheta'(\vartheta)+\zeta(\alpha)u+(\betam-\zeta(2\alpha)+1) u^2+o_\vartheta(1) \Big]_{\vartheta=\theta_u^*} =o_{\theta_u^*}(1).\nn
	}
Since $\vartheta=\thetastar_u$ is the solution of the variational problem in \refeq{thetastaru-def}, we have in particular,
$u^{\tau-1}\Fvartheta'(\thetastar_u)+\zeta(\alpha)u+(\betam-\zeta(2\alpha)+1) u^2=0$. Also recall Lemma \ref{lem-unique-theta} and Lemma \ref{lem-approx-VP} to see that $\theta_u^*$ is bounded away from zero for $u$ big enough. This shows that $u\widetilde \expec[\SS_{u}]=o(1)$, as required. \qed
\medskip

Recall the expression of $\SS_t$ in \eqref{SS-def-ref}. In order to investigate the
asymptotics of $\widetilde\expec[\SS_{t}]$, we start by describing the distribution
of the indicator processes $(\II_i(t))_{t\geq 0}$ under the measure $\widetilde \prob$.
Since our indicator processes $(\II_i(t))_{t\geq 0}$ are \emph{independent},
this property also holds under the measure $\widetilde \prob$:

\begin{Lemma}[Indicator processes under the tilted measure]
\label{lem-ind-tilt}
Under the measure $\widetilde \prob$, the distribution of the indicator processes
$(\II_i(t))_{t\geq 0}$ is that of independent indicator processes. More precisely,
    \eqn{
    \lbeq{II-proc-tilt}
    \II_i(t)=\indic{T_i\leq t},
    }
where $(T_i)_{i\geq 2}$ are independent random variables with distribution
    \eqn{
    \lbeq{Ti-def-tilt}
    \widetilde \prob(T_i\leq t) =
    \begin{cases}
    \frac{\e^{\theta c_iu}(1-\e^{-c_i t})}{\e^{\theta c_iu}(1-\e^{-c_iu})
    +\e^{-c_iu}} &\text{for } t \leq u;\\
    \frac{\e^{\theta c_iu}(1-\e^{-c_iu})
    +(\e^{-c_iu}-\e^{-c_i t})}{\e^{\theta c_iu}(1-\e^{-c_iu})
    +\e^{-c_iu}} &\text{for } t>u.
    \end{cases}
    }
\end{Lemma}
\bigskip

\noindent
{\it Proof of Lemma \ref{lem-expec-S}.} Part (d) of Lemma \ref{lem-expec-S} follows from Lemma \ref{extralemma}. It remains to prove Parts (a)--(c).

\noindent
Recall the definition of $I_{\sss E}(a)$ in \eqref{def_I_E}.
We calculate, by \eqref{Ti-def-tilt} and \eqref{SS-def-ref}, for $t \leq u$ with $\theta = \theta_u^*$,
    \eqn{
    \lbeq{expect}
    \widetilde \expec[\SS_t] - 1 -\betam t
    = \sum_{i=2}^{\infty} c_i(\widetilde \expec[\II_i(t)] - c_i t)
    = \sum_{i=2}^{\infty} c_i\Big( \frac{ \e^{\theta c_i u} ( 1 - \e^{-c_i t} ) }{ \e^{\theta c_i u}(1 - \e^{-c_i u}) + \e^{-c_i u} }  - c_i t\Big).
    }
For $y\ge 0$ and $a\in[0,1]$, define
\eqn{
\fff_1(y,a):= y\Big( \frac{ \e^{\theta y} ( 1 - \e^{-y a} ) }{ \e^{\theta y}(1 - \e^{- y}) + \e^{- y} }  - y a\Big).
}

\noindent  Equation (\ref{expect}) can be rewritten as
    \eqn{
    \lbeq{expect1}
    \widetilde \expec[\SS_t] - 1 -\betam t
    =
    {1\over u}\sum_{i=2}^{\infty} \fff_1\left(c_i u,t/u\right).
    }

\noindent We deduce also that for all $t \leq u$
    \eqn{
    \lbeq{expect2}
    \widetilde \expec[\SS_t-\SS_u] -\betam (t-u)
    =
    {1\over u}\sum_{i=2}^{\infty} \fff_2\left(c_i u,t/u\right)
    }

\noindent with $\fff_2(y,a):=\fff_{1}(y,a)-\fff_1(y,1)$. Moreover, we can write $\fff_1(y,a)=y(h(y,a) -ay)$ with $h(y,a)={h_1(ay) \over h_2(y)}$ and
	\eqn{
	h_1(z)= 1 - \e^{- z },
	\qquad
	h_2(z)= 1 - \e^{- z} + \e^{- (1+\theta)z}.
	}
Remember that $c_i=i^{-\alpha}$. To apply Lemma \ref{lem-taylor}, we need to control $\partial_y \fff_1(y,a)$ and $\partial_y \fff_2(y,a)$.
We have
\eqn{\partial_y h(y,a) = a{h_1'(ay) \over h_2(y)} - h_1(ay){h_2'(y) \over h_2^2(y)}.}
 Notice that $\Big|{h_1'(ay) \over h_2(y)}\Big| \le c  $ for any $y\ge 0$. Here we used that $\theta_u^*$ is uniformly bounded (see for instance Lemma \ref{lem-approx-VP}). On the other hand, $| h_1(ay){h_2'(y) \over h_2^2(y)}| \le c \e^{-y}$. This yields that
	\eqn{
	|\partial_y h(y,a) | \le c a + c\e^{-y}.
	}
Moreover, observe that $|h(y,a)| = |{h_1(ay) \over h_2(y)}| \le  cay$. Going back to $\fff_1(y,a)=y(h(y,a)-ay)$, we get
	\eqn{
	|\partial_y \fff_1(y,a)| \le y |\partial_y h(y,a)| + |h(y,a)| + 2ay \le
	cay + c.
	}
We further have
	\eqn{
	|\partial_y \fff_2(y,a)| \le|h(y,1)-h(y,a)| + y|\partial_y (h(y,1)-h(y,a))| + 2(1-a)y,
	}
so that $|h(y,1)-h(y,a)|=|{h_1(y)-h_1(ay)\over h_2(y)}| \le c(1-a)y$. Also,
$|\partial_y (h(y,1)-h(y,a))|\le |{h_1'(y) \over h_2(y)}| + a|{h_1'(ay) \over h_2(y)}| + |h_1(y)-h_1(ay)| |{h_2'(y) \over h_2^2(y)}|$ which is less than $c(\e^{-y} +a\e^{-ay} + \e^{-y})$. Therefore, for $a\in[1/2,1]$,
	\eqn{
	|\partial_y \fff_2(y,a)| \le c ( (1-a)y + y\e^{-y/3}).
	}
We can now use a straightforward extension of Lemma \ref{lem-taylor}
with $\gamma\in\{0,1\},b=0$ to see that there exists a constant $c>0$ such that, for any $u\ge 1$, and any $a\in [0,1]$,
	\eqn{
	\Big|\sum_{i=2}^{\infty} [ \fff_1(c_i u,a)-
	\int_i^{i+1} \fff_1(ux^{-\alpha},a)dx ] \Big| \le ca u^2 + c u.
	}
With $\gamma=1$ and $b \in \{0,1/3\}$ we obtain for any $u\ge 1$ and any $a\in[1/2,1]$
	\eqn{
	\Big|\sum_{i=2}^{\infty} [ \fff_2(c_i u,a)-
	\int_i^{i+1} \fff_2(ux^{-\alpha},a)dx ] \Big| \le c((1-a) u^2 + 1) .
	}
By (\ref{expect1}) and (\ref{expect2}), it follows that for $u\ge 1$ and $t\in[0,u]$,
	\eqn{
	\Big|\widetilde \expec[\SS_t] - 1 -\betam t
    	-
    	u^{-1}\int_{x\ge 2} \fff_1(ux^{-\alpha},t/u)dx\Big| \le c (t+1) ,
	}
and, for $u\ge 1$ and $t\in [u/2,u]$,
	\eqn{
	\Big|\widetilde \expec[\SS_t-\SS_u]  -\betam (t-u)
    	-
    	u^{-1}\int_{x\ge 2} \fff_2(ux^{-\alpha},t/u)dx\Big| \le c((u-t)  + u^{-1}) .
	}
Observe that
	\eqan{
	u^{-1}\int_{x\ge 2} \fff_1(ux^{-\alpha},t/u)dx
	&=
	(\tau-1)u^{\tau-2}\int_{0}^{2^{-\alpha}u} \fff_1(x,t/u)x^{-\tau}dx\\
	&=
	u^{\tau-2}\Big(I_{\sss E}^u(t/u) - (\tau-1)\int_{x\ge 2^{-\alpha}u} \fff_1(x,t/u)x^{-\tau}dx\Big), 	
	\nn
	}
where $I_{\sss E}^u$ is defined by replacing $\theta^*$ by $\theta_u^*$ in the definition of $I_{\sss E}$. We have as well
	\eqn{
	u^{-1}\int_{x\ge 2} \fff_2(ux^{-\alpha},t/u)dx
	=
	u^{\tau-2}(I_{\sss E}^u(t/u)-I_{\sss E}^u(1) - (\tau-1)\int_{x\ge 2^{-\alpha}u} \fff_2(x,t/u)x^{-\tau}dx).
	}
We have seen that $|\fff_1(x,a)|=x|h(x,a) -ax|\le c ax^2 $, which implies that $\int_{x\ge 2^{-\alpha}u} |\fff_1(x,t/u)|x^{-\tau}dx \le c tu^{2-\tau}$. Similarly, $|\fff_2(x,a)| \le c(1-a)x^2$ implies that $\int_{x\ge 2^{-\alpha}u} |\fff_2(x,t/u)|x^{-\tau}dx \le c (u-t)u^{2-\tau}=o(u-t)$. Consequently, for $u\ge 1$ and $t\in[0,u]$,
	\eqn{
	|\widetilde \expec[\SS_t]
    	-
   	u^{\tau-2} I_{\sss E}^u(t/u)| \le  c  (t+1),
	}
and for $u\ge 1$ and $t\in[u/2,u]$
	\eqn{
	\label{diff-means}
	|\widetilde \expec[\SS_t-\SS_u]
    	-
   	u^{\tau-2} (I_{\sss E}^u(t/u) - I_{\sss E}^u(1)) | \le  c ((u-t) + u^{-1}).
	}
Equation \eqref{diff-means} immediately allows us to prove Lemma \ref{lem-expec-S}(c). Indeed,
the function $a\mapsto I_{\sss E}^u(a)$ is differentiable, so that we can approximate
	\eqn{
	u^{\tau-2} (I_{\sss E}^u(t/u) - I_{\sss E}^u(1))=u^{\tau-3} (t-u)\frac{d}{d a} I_{\sss E}^u(a)
	\big|_{a=a^*},
	}
for some $a^*\in [t/u,1]$. Since $a^*$ is close to $1$ and $\thetastar_u$ is close to $\thetastar$,
$\frac{d}{d a} I_{\sss E}^u(a)\big|_{a=a^*}=I_{\sss E}'(1)+o(1)$.  Lemma \ref{lem-expec-S}(c)
follows once we note that also $(u-t)=o((u-t)u^{\tau-3})$.

To obtain Lemma \ref{lem-expec-S}(a) and (b), we apply (\ref{diff_I_E}) below and use $I_{\sss E}(1)=0$.
\qed
\bigskip

\noindent{\it Proof of Lemma \ref{lem-var-S-rep}.} 	
We similarly define $I_{\sss E}^u$, $I_{\sss V}^u$, $J_{\sss V}^u$ and $G_{\sss V}^u$ by replacing $\theta^*$ by $\theta_u^*$ in the definitions \eqref{def_I_E}-\eqref{GV-def}, and we check that, for any $a\in[0,1]$ and $u\ge 1$
	\eqan{
	\lbeq{diff_I_E}
    	|I_{\sss E}^u(a)-I_{\sss E}(a)| &\le c |\theta^* - \theta_u^*| a,\\
	\lbeq{diff_I_V}
	|I_{\sss V}(a)-I_{\sss V}^u(a)|   &\le c a |\theta^*-\theta_u^*|,\\
	\lbeq{diff_J_V}|J_{\sss V}(a) - J_{\sss V}^u(a)| &\le c (1-a) |\theta^* - \theta_u^*|,\\
	\lbeq{diff_G_V} |G_{\sss V}(a) - G_{\sss V}^u(a)| &\le c \min(a,1-a) |\theta^* - \theta_u^*|.
	}
To calculate the variance of $\SS_{t}$ under $\widetilde \prob$ for $t \in [0,u]$, recall that under $\widetilde \prob$ the indicator processes in the definition of $\SS_t$ in \eqref{SS-def-ref} are independent. We obtain, using \eqref{Ti-def-tilt},
    	\eqn{
    	\widetilde{{\rm Var}}[\SS_{t}]
    	= \sum_{i=2}^{\infty} c_i^2 \widetilde \expec[ (\II_i(t)-\widetilde \expec[\II_i(t)])^2 ]
    	= \sum_{i=2}^{\infty} c_i^2 \widetilde \prob(T_i \leq t) (1-\widetilde \prob(T_i \leq t))
    	= u^{-2} \sum_{i=2}^{\infty} \fff_3\left(c_i u,t/u\right),  \lbeq{sum-for-var}
    	}
where
	\eqn{
	\fff_3(y,a):=
	y^2 \frac{ \e^{\theta y} ( 1 - \e^{- ay} ) }{ \e^{\theta y}(1 - \e^{- y}) + \e^{- y} }
	\left( 1 - \frac{ \e^{\theta y} ( 1 - \e^{- a y} ) }{ \e^{\theta y}(1 - \e^{-y}) + \e^{-y} } \right)
	}
with $\theta = \theta_u^*$. We have as well
	\eqn{
    	\widetilde{{\rm Var}}[\SS_{t}-\SS_{u}]
    	= \sum_{i=2}^{\infty} c_i^2 \widetilde \prob(T_i \in (t,u]) (1-\widetilde \prob(T_i \in (t,u]))
    	= u^{-2} \sum_{i=2}^{\infty} \fff_4\left(c_iu,t/u\right) \lbeq{sum-for-var-2}
    	}
with
    	\eqn{
    	\fff_4(y,a)= y^2 \frac{ \e^{\theta y} ( \e^{- ay} - \e^{-y}) }{ \e^{\theta y}(1 - \e^{- y}) + \e^{- y} }
	\left( 1 - \frac{ \e^{\theta y} (  \e^{- a y} -\e^{-y}) }{ \e^{\theta y}(1 - \e^{-y}) + \e^{-y} } \right)
    	}
and
	\eqn{\lbeq{sum-for-var-3}
	\widetilde{ {\rm Cov}}[\SS_t,\SS_u-\SS_t] = - \sum_{i=2}^\infty c_i^2
	\widetilde \prob(T_i\le t)\widetilde \prob(T_i\in(t,u]) = -u^{-2} \sum_{i=2}^\infty \fff_5\left(c_iu,t/u\right)
	}
with
	\eqn{
	\fff_5(y,a)
	= y^2{\e^{2\theta y} (1-\e^{-ay})(\e^{-ay} -\e^{-y}) \over (\e^{\theta y}(1-\e^{-y}) + \e^{-y})^2}.
	}
Let again $h(y,a):={h_1(ay)\over h_2(y)}$ where $h_1(z):=1-\e^{-z}$ and $h_2(z):=1-\e^{-z}+\e^{-(1+\theta)z}$. Then,
	\eqan{
	\fff_3(y,a) & = y^2h(y,a)(1-h(y,a)),\\
	\fff_4(y,a) & = y^2(h(y,1)-h(y,a)) (1- h(y,1)+h(y,a)),\\
	\fff_5(y,a) & = y^2h(y,a)(h(y,1)-h(y,a)).
	}
We bound $\partial_y \fff_i(y,a)$ for $i\in\{3,4,5\}$. We have
	\eqn{
	\partial_y h(y,a) = a{h_1'(ay) \over h_2(y)} - h_1(ay){h_2'(y) \over h_2^2(y)}.
	}
Firstly, $h_1'(ay)=\e^{-ay}$. Since $h_2\ge c>0$ (recall that $\theta_u^*$ is uniformly bounded), we get $0 \leq h(y,a) \leq c$ and $|a{h_1'(ay) \over h_2(y)}|\le ca\e^{-ay}$. Secondly, $|h_1(ay)|\le 1$ and $|h_2'(y)|\le (2+\theta)\e^{-y}$ hence
$|h_1(ay){h_2'(y) \over h_2^2(y)}| \le c\e^{-y}$. We get that
	\eqn{
	\lbeq{bound-dy-h}
	|\partial_y h(y,a)| \le ca\e^{-ay} + c\e^{-y}.
	}
From  $\partial_y \fff_3(y,a) = 2yh(y,a)(1-h(y,a)) +y^2\partial_y h(y,a)(1-2h(y,a))$, it now follows that
	\eqn{
	|\partial_y \fff_3(y,a)| \le c(1+y).
	}
Similarly,
	\eqn{
	\partial_y (h(y,a)-h(y,1)) = {ah_1'(ay) -h_1'(y)\over h_2(y)} - (h_1(ay)-h_1(y)){h_2'(y) \over h_2^2(y)}.
	}
We have $|h_1(ay)-h_1(y)|\le (1-a)y\e^{-ay}$, $|{h_2'(y) \over h_2^2(y)}|\le c$ and $|ah_1'(ay) -h_1'(y)|\le c(1-a)(1+y)\e^{-ay}$. This gives
	\eqn{
	|\partial_y (h(y,a)-h(y,1))| \le c (1-a)(1+y)\e^{-ay}.
	}
Since
	\eqn{
	\lbeq{bound-dy-f4}
	|\partial_y \fff_4(y,a)|
	\le
	c(y|h(y,a)-h(y,1)| + y^2|\partial_y (h(y,a)-h(y,1))|),
	}
we deduce that, for $a\in[1/2,1]$
	\eqn{
	|\partial_y \fff_4(y,a)|
	\le
	c(1+y^3)(1-a)\e^{-ay} \le c(1-a).
	}
On the other hand, if $a\in[0,1/2]$, we write, this time using \eqref{bound-dy-h} to bound the second term in \eqref{bound-dy-f4},
	\eqn{
	|\partial_y \fff_4(y,a)|
	\le
	c(1+y).
	}
Similarly,
	\eqan{
	|\partial_y \fff_5(y,a)| \le 2y|h(y,a)(h(y,1)-h(y,a))| + y|y \partial_y h(y,a)||h(y,1)-h(y,a)|
	+ y^2 h(y,a)|\partial_y (h(y,1)-h(y,a))|.
	}
We use that the terms $|h(y,a)|, |h(y,1)-h(y,a)|, |y\partial_y h(y,a)|$ are bounded by a constant, and $y^2|\partial_y (h(y,a)-h(y,1))| \le c(1-a)$ if $a\in[1/2,1]$ and $y^2|\partial_y (h(y,a)-h(y,1))| \le c(1+y)$ if $a\in[0,1/2]$. We get that
	\eqn{
	|\partial_y \fff_5(y,a)| \le c(1+y)\indic{a\in[0,1/2]} + c (1-a) (1+y) \indic{a\in[1/2,1]}.
	}
Next, we use Lemma \ref{lem-taylor} as before to see that
	\eqan{
	\Big|\sum_{i=2}^{\infty} [ \fff_3(u c_i,a)-\int_i^{i+1} \fff_3(u x^{-\alpha},a)dx ] \Big| &\le c (u+u^2),\\
	\Big|\sum_{i=2}^{\infty} [ \fff_4(u c_i,a)-\int_i^{i+1} \fff_4(u x^{-\alpha},a)dx ] \Big|
	&\le c(u+u^2) 	\indic{a\in[0,1/2]} + c(1-a)u \indic{a\in [1/2,1]},\\
	\Big|\sum_{i=2}^{\infty} [ \fff_5(u c_i,a)-\int_i^{i+1} \fff_5(u x^{-\alpha},a)dx ] \Big|
	&\le c(u+u^2) \indic{a\in[0,1/2]} + c(1-a)(u+u^2) \indic{a\in[1/2,1]}.
	}
Going back respectively to (\ref{sum-for-var}), (\ref{sum-for-var-2}) and (\ref{sum-for-var-3}), we get
	\eqan{
	\Big| \widetilde{{\rm Var}}[\SS_{t}] - u^{-2} \int_{2}^\infty \fff_3(u x^{-\alpha},t/u)dx\Big| &\le c, \\
	\Big| \widetilde{{\rm Var}}[\SS_{t}-\SS_u] - u^{-2} \int_{2}^\infty \fff_4(u x^{-\alpha},t/u)dx\Big| &\le
	c \indic{t/u\in[0,1/2]} + c(u-t) u^{-2} \indic{t/u\in [1/2,1]} \le c (1-t/u), \\
	\Big| \widetilde{{\rm Cov}}[\SS_{t},\SS_u-\SS_t] + u^{-2} \int_{2}^\infty \fff_5(u x^{-\alpha},t/u)dx\Big|
	 &\le
	c \indic{t/u\in[0,1/2]} + c(u-t)u^{-1} \indic{t/u\in[1/2,1]} \le c (1-t/u).
	}
With the change of variables $y=x^{-\alpha}u$ we see that (recall that $\theta=\theta_u^*$)
	\eqan{
	\int_{2}^\infty \fff_3(u x^{-\alpha},t/u)dx
	& =
	u^{\tau-1}\Big(I_{\sss V}^u(t/u) - (\tau-1) \int_{y\ge 2^{-\alpha}u} y^{-\tau}\fff_3(y,t/u)dy\Big),\\
	\int_{2}^\infty \fff_4(u x^{-\alpha},t/u)dx
	&=
	u^{\tau-1}\Big(J_{\sss V}^u(t/u) - (\tau-1) \int_{y\ge 2^{-\alpha}u} y^{-\tau}\fff_4(y,t/u)dy\Big),\\
	\int_{2}^\infty \fff_5(u x^{-\alpha},t/u)dx
	&=
	u^{\tau-1}\Big(G_{\sss V}^u(t/u) - (\tau-1) \int_{y\ge 2^{-\alpha}u} y^{-\tau}\fff_5(y,t/u)dy\Big).
	}
Since $|\fff_3(y,a)|\le cy^2$, we have $|\int_{y\ge 2^{-\alpha}u} y^{-\tau}\fff_3(y,t/u)dy| \le u^{3-\tau}$. We arrive at
	\eqn{
	\Big| \widetilde{{\rm Var}}[\SS_{t}] - u^{\tau-3}I_{\sss V}^u(t/u) \Big| \le c.
	}
We check similarly that $|\fff_4(y,a)|\le c (1-a) y^2$ and $|\fff_5(y,a)| \le c(1-a)y^2$, so $|\int_{y\ge 2^{-\alpha}u} y^{-\tau}f(y,t/u)dy| \le c (u-t)u^{2-\tau}$ for $\fff$ being $\fff_4$ or $\fff_5$. Therefore, for $t\in [0,u]$,
	\eqan{
	\Big| \widetilde{{\rm Var}}[\SS_{t}-\SS_u] - u^{\tau-3}J_{\sss V}^u(t/u)\Big| &\le c (u-t) u^{-1},\\
	\Big| \widetilde{{\rm Cov}}[\SS_{t},\SS_{u}-\SS_t] - u^{\tau-3}G_{\sss V}^u(t/u)\Big| &\le c(u-t)u^{-1}.
	}
Use (\ref{diff_I_V}), (\ref{diff_J_V}) or (\ref{diff_G_V}) to complete the proof.
\hfill\qed
\medskip

We next investigate what happens to the means and variances for small $a$ or for $a$ close to 1:
\begin{Lemma}[Asymptotic mean and variance near extremities]
\label{lem-mean-var-extreme}
\ \\
(a) For $a \in [0,1]$, $I_{\sss E}(a) \leq a I_{\sss E}'(0)$ and
$I_{\sss E}(a) \leq -(1-a)I_{\sss E}'(1)$.\\
(b) As $a\to 1$,  $J_{\sss V}(a) = -(1-a)J_{\sss V}'(1)(1+o(1))$ with $J_{\sss V}'(1)<0$, while, as $a\to 0$,
	\eqn{
	\lbeq{I_V-a-near-zero}
	I_{\sss V}(a) = a^{\tau-3} \mathcal I_{\sss V} (1+o(1)),
	\qquad
	\mbox{with}
	\qquad
	\mathcal I_{\sss V}= (\tau-1) \int_0^\infty (  1 - \e^{-y}  ) \e^{-y} \frac{ dy }{ y^{\tau - 2} }.
	}
Consequently, there exist $0<\underline{c}<\overline{c}<\infty$ such that, for every $a\in[0,\vep]$ with $\vep>0$ sufficiently small,
	\eqn{
	\lbeq{I_V-a-near-zero-rep}
	 \underline{c} a^{\tau-3} \leq I_{\sss V}(a) \leq \overline{c} a^{\tau-3}.
	}
\end{Lemma}

\proof
For the proof of Part (a), use the concavity of $I_{\sss E}(a)$ on $[0,1]$. For Part (b), note that $J_{\sss V}$ is continuously differentiable on $(0,1]$. Equation \eqref{I_V-a-near-zero} follows by the change of variable $x:= av$ in \eqref{IV-def} and dominated convergence. Equation \refeq{I_V-a-near-zero-rep} follows directly from
\refeq{I_V-a-near-zero} and the fact that $a\in[0,\vep]$ with $\vep>0$ sufficiently small.
\qed
\medskip

\section{Joint distribution of $\SS_t$ and $\SS_u$: proof of Proposition \ref{prop-laplace}}
\label{profff}
The proof follows by explicitly computing the joint moment generating function
of $(\SS_t,\SS_u-\SS_t)$ using Lemma \ref{lem-ind-tilt} and studying its asymptotics.
We prove parts (a) and (b) simultaneously, by noting that the extra assumption $t\not\in (u-u^{-(\tau-5/2)},u]$ is not needed when $\lambda_2=0$.

We start by introducing some notation.
Let $\lambda_1,\lambda_2$ be elements of compact sets, and abbreviate (recall that we consider $t\in[\vep,u-u^{-(\tau-5/2)}]$)
	\eqan{
	\lbeq{lambda_12}
	\tilde \lambda_1 &= \lambda_1 I_{\sss V}(t/u)^{-1/2} u^{-(\tau-3)/2},
	\qquad\quad
	\tilde \lambda_2 = \lambda_2 J_{\sss V}(t/u)^{-1/2} u^{-(\tau-3)/2}.
	}
For Part (a), we simply take $\tilde \lambda_2\equiv0$. For these choices, \eqref{IV-domain}, \eqref{JV-domain} and Lemma \ref{lem-mean-var-extreme}(b) guarantees that there exists a $c>0$ such that
	\eqn{
	\lbeq{tilde-lambda12-bds}
	|\tilde\lambda_1| \le c t^{-(\tau-3)/2},
	\qquad
	|\tilde \lambda_2| \le c \sqrt{\frac{u}{u-t}} u^{-(\tau-3)/2}.
	}
We observe that
	\eqn{
	\widetilde \expec\left[\e^{\tilde \lambda_1 (\SS_{t} - \widetilde \expec[\SS_t])
	+  \tilde \lambda_2 (\SS_{u} - \SS_t - \widetilde \expec[\SS_t - \SS_u])}\right]
	=
	\prod_{j\ge 2} \widetilde \expec\left[\e^{ c_j ( \tilde \lambda_1 (\II_j(t) - \widetilde \expec[\II_j(t)])
	+ \tilde \lambda_2 (\II_j(u)-\II_j(t) - \widetilde \expec[\II_j(u)-\II_j(t)]) )}  \right].
	}
From the distribution of the indicators $(\II_j(t))_{j\ge 2}$ under $\widetilde \prob$ given in Lemma \ref{lem-ind-tilt}, we get that
	\eqn{\label{eq:fourier}
	\log \widetilde \expec\left[\e^{\tilde \lambda_1 (\SS_{t} - \widetilde \expec[\SS_t])
	+  \tilde \lambda_2 (\SS_{u} - \SS_t - \widetilde \expec[\SS_t - \SS_u])}\right]
	= \sum_{j\ge 2} a_j,}
where we define
	\eqn{\label{def:aj}
	a_j =\log \left( 1+{(\e^{\tilde \lambda_1 c_j}-1) (1-\e^{-t c_j})
	+ (\e^{\tilde \lambda_2 c_j}-1) (\e^{-c_j t} - \e^{-c_j u})\over 1 - \e^{-uc_j}
	+\e^{-u c_j (1+\theta)} } \right) - c_j {\tilde \lambda_1(1-\e^{-c_jt})
	+ \tilde \lambda_2(\e^{-c_j t}-\e^{-c_j u}) \over 1-\e^{-uc_j} + \e^{-uc_j(1+\theta)}},
	}
with $\theta=\theta_u^*$. A Taylor expansion of $\log(1+x)$ around~$1$ (we use that $|\log(1+x)-x+{x^2\over 2}| \le c |x|^3$ for all $x$ greater than some $\rho>-1$) shows that
	\eqn{
	a_j = (\e^{\tilde \lambda_1 c_j}-1-c_j\tilde \lambda_1) {1-\e^{-t c_j} \over  1 - \e^{-uc_j}
	+\e^{-u c_j 	(1+\theta)}} + (\e^{\tilde \lambda_2 c_j}-1-\tilde \lambda_2 c_j)
	{ \e^{-c_j t} - \e^{-c_j u} \over 1 - \e^{-uc_j} +\e^{-u c_j (1+\theta)} }  - {1\over 2}b_j^2 + O(|b_j|^3),
	}
where we let
	\eqn{
	b_j = {(\e^{\tilde \lambda_1 c_j}-1) (1-\e^{-t c_j}) + (\e^{\tilde \lambda_2 c_j}-1)
	(\e^{-c_j t} - \e^{-c_j u})\over 1 - \e^{-uc_j} +\e^{-u c_j (1+\theta)} }.
	}
Here, we claim that $b_j$ is uniformly greater than $-1$. Indeed,  for $t\ge  \vep$,
\refeq{tilde-lambda12-bds} implies that $|\tilde \lambda_1| \le c(\vep)$,
while $|\tilde \lambda_2|\le c ( (u-t)u^{\tau-4})^{-1/2}$.
Hence, for $t \leq (1-\epsilon) u$, $|\tilde \lambda_2| = o_u(1)$ and
$|\e^{\tilde \lambda_2 c_j}-1| (\e^{-c_j t} - \e^{-c_j u}) = o_u(1)$. If $t > (1-\vep) u$, it remains to investigate the case where $\tilde \lambda_2 < 0$. This case is absent for part (a), since $\tilde \lambda_2= 0$. We obtain for $t \leq u-u^{-(\tau-5/2)}$,
	\eqn{
	|\e^{\tilde \lambda_2 c_j}-1| (\e^{-c_j t} - \e^{-c_j u})
	\leq c|\lambda_2| u^{-1/4} (u c_j) \e^{-c_j (1-\vep) u}
	= o_u(1).
	}
This completes the proof that $b_j > -1$.

Using that for $i=1,2$, $\e^{\tilde \lambda_i c_j}-1-c_j\tilde \lambda_i = {1\over 2}c_j^2 \tilde \lambda_i^2 + O(c_j^3|\tilde \lambda_i|^3)$, we get
	\eqan{
	\lbeq{taylor-for-aj}
	a_j =& {c_j^2 \over 2} (\tilde \lambda_1)^2 {1-\e^{-t c_j} \over
	1 - \e^{-uc_j} + \e^{-u c_j (1+\theta)}} + {c_j^2 \over 2} (\tilde \lambda_2)^2 {\e^{-t c_j}-\e^{-u c_j}
	 \over  1 - \e^{-uc_j} +\e^{-u c_j (1+\theta)}} - {1\over 2}b_j^2 + O(|b_j|^3) \\
  	& \quad+ O\Big( |c_j\tilde \lambda_1|^3+|c_j \tilde \lambda_2|^3 (\e^{-c_j t}-\e^{-c_ju}) \Big). \nn
	}

We continue to investigate the term $b_j^2$. We develop $b_j^2$ so that
	\eqn{
	b_j^2 = (\e^{\tilde \lambda_1 c_j} -1)^2 b_{j,1} + (\e^{\tilde \lambda_2 c_j} -1)^2 b_{j,2} +
	(\e^{\tilde \lambda_1 c_j}-1)(\e^{\tilde \lambda_2 c_j} -1)b_{j,3},
	}
with the obvious notation. Our aim is to control the error
$b_j^2 - \big((\tilde \lambda_1 c_j)^2 b_{j,1} + (\tilde \lambda_2 c_j)^2 b_{j,2} +
	\tilde \lambda_1 \tilde \lambda_2 c_j^2 b_{j,3}\big).$
For this, we use the approximations $(\e^{\tilde \lambda_i c_j} -1)^2 = (c_j \tilde \lambda_i)^2 + O((c_j\tilde \lambda_i)^3)$ for $i=1,2$, as well as the fact that $\inf_{u\geq 0}[1 - \e^{-uc_j} +\e^{-u c_j (1+\theta)}]>0$, to see that
	\eqan{
	& [(\e^{\tilde \lambda_1 c_j} -1)^2 -(\tilde \lambda_1 c_j)^2]b_{j,1}
	+ [(\e^{\tilde \lambda_2 c_j} -1)^2 - (\tilde \lambda_2 c_j)^2]b_{j,2} \\
	&=  O\Big( (c_j\tilde \lambda_1)^3+(c_j \tilde \lambda_2)^3 (\e^{-c_j t}-\e^{-c_ju})^2 \Big). \nn
	}
We next bound $[(\e^{\tilde \lambda_1 c_j}-1)(\e^{\tilde \lambda_2 c_j} -1)-\tilde \lambda_1 \tilde \lambda_2 c_j^2 ]b_{j,3}$. We use that for $i=1,2$, $(\e^{\tilde \lambda_i c_j}-1) = c_j \tilde \lambda_i + O((c_j \tilde \lambda_i)^2)$. As before,
	\eqan{
	& [(\e^{\tilde \lambda_1 c_j}-1)(\e^{\tilde \lambda_2 c_j} -1)
	-\tilde \lambda_1 \tilde \lambda_2 c_j^2 ]b_{j,3}\\
	&=
	O\Big( \Big[ c_j^3 \Big\{ |\tilde\lambda_1| (\tilde \lambda_2)^2 + (\tilde \lambda_1)^2 |
	\tilde\lambda_2| \Big\} + c_j^4 (\tilde \lambda_1)^2 (\tilde \lambda_2)^2 \Big]
	(\e^{-c_j t}-\e^{-c_ju}) \Big). \nn
	}
Therefore, we finally arrive at
	\eqan{
	\lbeq{order-term-error-bj2}
	& b_j^2 - \Big((\tilde \lambda_1 c_j)^2 b_{j,1} + (\tilde \lambda_2 c_j)^2 b_{j,2}
	+ \tilde \lambda_1 \tilde \lambda_2 c_j^2 b_{j,3}\Big) \\
	&= c_j^3 O\Big((\tilde \lambda_1)^3+ (\e^{-c_j t}-\e^{-c_ju})
	\Big\{(\e^{-c_j t}-\e^{-c_ju}) (\tilde \lambda_2)^3
	+ |\tilde\lambda_1| (\tilde \lambda_2)^2 + (\tilde \lambda_1)^2 |\tilde\lambda_2|
	+ c_j (\tilde \lambda_1)^2 (\tilde \lambda_2)^2 \Big\}\Big). \nn
	}

A similar reasoning as the above shows, using once again that $|\tilde \lambda_1| \le c(\vep)$, that the remaining term $O(|b_j|^3)$ can be bounded by
	\eqn{
	\lbeq{order-term-error-bj3}
	O(|b_j|^3) \leq
	O\Big( \sum_{k=0}^3 |c_j \tilde \lambda_1|^{3-k} |c_j \tilde \lambda_2|^k
	(\e^{-c_j t}-\e^{-c_ju})^k \Big).
	}
Therefore we can write
	\eqn{
	a_j = c_j^2(\tilde \lambda_1^2 a_{j,1} + \tilde \lambda_2^2 a_{j,2} -
	\tilde \lambda_1\tilde \lambda_2 a_{j,3}) + \Theta_j,
	}
where
	\eqan{
	a_{j,1} &=  {1\over 2}{ 1-\e^{-t c_j}\over 1 - \e^{-uc_j} +\e^{-u c_j (1+\theta)} }
	\left(1 -{ 1-\e^{-t c_j}\over 1 - \e^{-uc_j} +\e^{-u c_j (1+\theta)} }\right),\\
	a_{j,2} &= {1\over 2}{ \e^{-t c_j} -\e^{-u c_j} \over 1 - \e^{-uc_j} +\e^{-u c_j (1+\theta)} }
	\left(1- { \e^{-t c_j} -\e^{-u c_j} \over 1 - \e^{-uc_j} +\e^{-u c_j (1+\theta)} }\right),\\
	a_{j,3} &= { (1-\e^{-t c_j})(\e^{-t c_j} -\e^{-u c_j})
	\over (1 - \e^{-uc_j} +\e^{-u c_j (1+\theta)})^2 },
	}
and, by collecting terms from \eqref{taylor-for-aj}, \eqref{order-term-error-bj2} and \eqref{order-term-error-bj3}, the error term $\Theta_j$ satisfies
	\eqn{
	\label{Theta-j-bd}
	|\Theta_j| \le  c\,c_j^3 \Big( |\tilde \lambda_1|^3
	+ (\e^{-c_j t}-\e^{-c_ju}) \Big[ |\tilde \lambda_2|^3
	+ |\tilde \lambda_1| (\tilde \lambda_2)^2
	+ (\tilde \lambda_1)^2 |\tilde \lambda_2| + c_j (\tilde \lambda_1)^2 (\tilde \lambda_2)^2 \Big] \Big).
	}

We continue to bound $\sum_{j\geq 2}|\Theta_j|$.
We use \refeq{tilde-lambda12-bds} to obtain
	\eqn{
	\sum_{j \geq 2} c_j^3 |\tilde \lambda_1|^3 \leq c t^{-3(\tau-3)/2}.
	}
This completes the bound on $\sum_{j\geq 2}|\Theta_j|$ for part (a). For part (b) and the other terms, we split, depending on whether $t \leq (1-\vep) u$ or not. In the case $t \leq (1-\vep) u$ for arbitrary $\vep>0$, $|\tilde \lambda_2| \leq c(\vep) u^{-(\tau-3)/2} = o_u(1)$ and thus
	\eqn{
	\sum_{j \geq 2} |\Theta_j| \leq c t^{-3(\tau-3)/2} + o_u(1).
	}
In the case $t > (1-\vep) u$ we use
the bound on $|\tilde \lambda_2|$ in \refeq{tilde-lambda12-bds}
and $t \leq u-u^{-(\tau-5/2)}$ to bound
	\eqan{
	\sum_{j \geq 2} c_j^3 |\tilde \lambda_2|^3 (\e^{-c_j t}-\e^{-c_ju})
	& \leq c \sum_{j \geq 2} c_j^3 u^{3/2} (u-t)^{-1/2} u^{-3(\tau-3)/2} \e^{-c_j t} c_j \\
	& \leq c \sum_{j \geq 2} c_j^3 u^{3/2} u^{(\tau-5/2)/2} u^{-3(\tau-3)/2}
	\e^{-c_j (1-\vep) u} c_j 	\nn\\
	& \leq c \sum_{j \geq 2} (u c_j)^4  u^{1-\tau} \e^{-c_j (1-\vep) u} u^{-1/4}= O(u^{-1/4}). \nn
	}
For the remaining terms when $t>(1-\vep)u$, we use $\e^{-c_j t}-\e^{-c_ju}\leq \e^{-c_j (1-\vep) u}
c_j (u-t)$ to finally obtain
	\eqan{
	\label{Theta-j-bd-end}
	\sum_{j \geq 2} |\Theta_j|
	&\le c t^{-3(\tau-3)/2} + o_u(1) +
	c\, \sum_{j \geq 2} c_j^3 \Big(\e^{-c_j (1-\vep) u}
	c_j (u-t) \Big[ (\tilde \lambda_2)^2
	+ |\tilde \lambda_2| \Big] \Big) \\
	&\le c t^{-3(\tau-3)/2} + o_u(1) + c\, \sum_{j \geq 2} c_j^3
	\Big(\e^{-c_j (1-\vep) u}
	\Big[(c_j u) u^{-(\tau-3)} + (c_j u) u^{-(\tau-3)/2}\Big)\Big]\nn\\
	&\le c t^{-3(\tau-3)/2} + o_u(1). \nn
	}

It follows that
	\eqn{
	\lbeq{log-mom-gen-funct}
	\log \widetilde \expec\left[\e^{\tilde \lambda_1 (\SS_{t} - \widetilde \expec[\SS_t])
	+ \tilde \lambda_2 (\SS_{u} - \SS_t - \widetilde \expec[\SS_t - \SS_u])}\right]
	= \sum_{j\ge 2} c_j^2(\tilde \lambda_1^2 a_{j,1} + \tilde \lambda_2^2 a_{j,2} -
	\tilde \lambda_1\tilde \lambda_2 a_{j,3}) +\Theta',
	}
where $\Theta'=\sum_{j \geq 2} \Theta_j$ satisfies $|\Theta'|\le c t^{-3(\tau-3)/2}+o_u(1)$.
We now estimate the sums over $j$ in
\refeq{log-mom-gen-funct}. We notice that (recall \eqref{sum-for-var},
\eqref{sum-for-var-2} and \eqref{sum-for-var-3})
	\eqn{
	\sum_{j\ge 2} c_j^2 a_{j,1} = {1\over 2}\widetilde {\rm Var}[\SS_t],\qquad
	\sum_{j\ge 2} c_j^2 a_{j,2} = {1\over 2}\widetilde {\rm Var}[\SS_t-\SS_u],\qquad
	\sum_{j\ge 2} -c_j^2 a_{j,3} = \widetilde {\rm Cov}[\SS_t,\SS_u-\SS_t].
	}
Lemmas \ref{lem-var-S-rep}, \ref{lem-mean-var-extreme}
and \ref{lem-approx-VP} imply that there exist $R_i=O(1)$ for $i\in \{1,2,3\}$
such that
	\eqan{
	\sum_{j\ge 2} c_j^2 a_{j,1} &= \frac{1}{2} I_{\sss V}(t/u)u^{\tau-3} + R_1,\\
	\sum_{j\ge 2} c_j^2 a_{j,2} &= \frac{1}{2} J_{\sss V}(t/u)u^{\tau-3} + \frac{u-t}{u} R_2,\\
	\sum_{j\ge 2} -c_j^2 a_{j,3} &= -G_{\sss V}(t/u)u^{\tau-3} + \frac{u-t}{u} R_3.
	}
We obtain, using \eqref{lambda_12}
	\eqn{
	\label{log-mom-gen-function}
	\log \widetilde \expec\left[\e^{\tilde \lambda_1 (\SS_{t} - \widetilde \expec[\SS_t])
	+  \tilde \lambda_2 (\SS_{u} - \SS_t - \widetilde \expec[\SS_t - \SS_u])}\right]
	=  {1\over 2}\lambda_1^2  +  {1\over 2}\lambda_2^2
	- \lambda_1\lambda_2 {G_{\sss V}(t/u) \over \sqrt{I_{\sss V}(t/u)J_{\sss V}(t/u)}} +\Theta,
	}
with
	\eqn{
	\Theta=\Theta'+\tilde\lambda_1^2 R_1+\tilde\lambda_2^2 (u-t)R_2/u
	+\tilde\lambda_1\tilde\lambda_2 (u-t)R_3/u.
	}
Using \refeq{tilde-lambda12-bds}, we see that the terms involving $R_i$ for $i\in\{1,2,3\}$ are
bounded by $ct^{-(\tau-3)}+cu^{-(\tau-3)}$,
so that
$|\Theta|\le o_u(1) + ct^{3(3-\tau)/2}$, which completes the proof.
\qed

\section{Density of $\SS_u$: proof of Proposition \ref{lem-dens-f} }
\label{sec-dens-Su}
In this section, we derive the asymptotics of the density $\widetilde f_{\SS_t}$ in Proposition \ref{lem-dens-f}.
We use the Fourier inversion formula
    \eqn{
    \lbeq{fourier}
    \widetilde f_{\SS_t}(s) =  \int_{-\infty}^{\infty} \widetilde \expec[\e^{2\ii \pi k  \SS_t}]\e^{-2\ii \pi k s}dk,
    }
where $\ii$ denotes the imaginary unit.
By a change of variables, we get
    \eqn{\label{eq:density_inverse}
     \widetilde f_{\SS_t}(s) =  u^{-(\tau-3)/2}\int_{-\infty}^{\infty} \widetilde \expec[\e^{2\ii \pi k  u^{-(\tau-3)/2}\SS_t}]\e^{-2\ii \pi k u^{-(\tau-3)/2} s}dk.
    }

We need to prove asymptotics when $t=u$ and an upper bound uniformly in $t\in[u/2,u]$. We will do both at the same time, and start by setting the stage for $t=u$. Remember that $\widetilde \E[\SS_u]=o(u^{-1})$ by our choice of $\theta$ (see Lemma \ref{extralemma}).
By Proposition \ref{prop-laplace}(a),
	\eqn{
	\label{conv-FT-Su}
	\widetilde \expec[\e^{2\ii \pi k  u^{-(\tau-3)/2}\SS_u}]\rightarrow \e^{-(2\pi k)^2I_{\sss V}(1)/2}.
	}

We want to use dominated convergence. The bound used for dominated convergence will then immediately prove the uniform upper bound for all $t\in[u/2,u]$. Write
	\eqn{
	r_t(k):= \Big|\widetilde \expec[\e^{2\ii \pi k  u^{-(\tau-3)/2}\SS_t}]\Big|
	=\Big|\widetilde \expec[\e^{2\ii \pi k  u^{-(\tau-3)/2}(\SS_t-\widetilde \expec[\SS_t])}]\Big|
	}
for the modulus of $\widetilde \expec[\e^{2\ii \pi k  u^{-(\tau-3)/2}(\SS_t-\widetilde \expec[\SS_t])}]$ for $t\in[u/2,u]$. We can compute $r_t(k)$ explicitly, using \eqref{Ti-def-tilt}. We find that
 	\eqn{
 	\log r_t(k)
	= \frac{1}{2} \sum_{j\ge 2}\log \Big\{ 1- 2 \e^{-c_jt}\e^{\theta c_ju}
	\frac{ (1-\e^{-c_ju}) (1-\cos (2\pi k c_j u^{-(\tau-3)/2})) }
	{ [\e^{-c_ju} + \e^{\theta c_ju}(1-\e^{-c_ju})]^2 } \Big\}.
	}
Using the inequality $\log(1-x)\le -x$ for $x<1$, it follows that
	\eqan{
    \lbeq{density-calc}
	\log r_t(k) &\le -\sum_{j\ge 2} \e^{-c_ju}\e^{\theta c_ju} \frac{ (1-\e^{-c_jt}) (1-\cos (2\pi k c_j u^{-(\tau-3)/2})) }{ [\e^{-c_ju} + \e^{\theta c_ju}(1-\e^{-c_ju})]^2 }\\
	&\le   - c\sum_{j\ge 2:c_j<\frac{1}{u}}  (1-\e^{-c_jt}) (1-\cos (2\pi k c_j u^{-(\tau-3)/2})). \nn
	}
We have $1-\e^{-x}\ge c x$ for $x\in[0,1]$. Hence (remember that $t\in[u/2,u]$)
	\eqn{
	\log r_t(k) \le -c u \sum_{j\ge 2:c_j< \frac{1}{u}} c_j  (1-\cos (2\pi k c_j u^{-(\tau-3)/2})) .
	}
We split the integral depending on the value of $k$. First suppose that $ k u^{-(\tau-1)/2}\le 1/8$ and $c$ such that $1-\cos(\pi x) \ge c x^2$ for any $x\in[0,\pi/4]$. Then $1-\cos (2\pi k c_j u^{-(\tau-3)/2}) \ge c k^2 c_j^2 u^{-(\tau-3)}$ for any $c_j <1/u$, so that
	\eqn{
    	\lbeq{density-1}
	\log r_t(k) \le -c u^{4-\tau} k^2 \sum_{j\ge 2:c_j< \frac{1}{u}} c_j^3 \le  - c k^2.
	}
To obtain the last inequality we used that $c_j\leq 1/u$ precisely when $j \geq u^{\tau-1}$ and that
	\eqan{
	u^{4-\tau} \sum_{j\colon c_j\leq 1/u} c_j^3
	&\geq u^{4-\tau} \int_{u^{\tau-1}}^\infty x^{-3/(\tau-1)} dx
	= u^{4-\tau} \left( -\tfrac{3}{\tau-1}+1 \right)^{-1} \left[ x^{-3/(\tau-1)+1} \right]_{u^{\tau-1}}^\infty \nn\\
	&= C(\tau) u^{4-\tau} u^{(\tau-1)(-3/(\tau-1)+1)}= C(\tau), \nn
	}
where we have used that $-3/(\tau-1)+1 \in (-1/2,0)$.
In the other case, let $y_{k}:= 8 k u^{-(\tau-3)/2}$ which is greater than $u$ by assumption. We have similarly
	\eqan{
    	\lbeq{density-2}
	\log r_t(k) &\le - cu \sum_{j\ge 2:c_j<\frac{1}{y_k}}  c_j  (1-\cos (2\pi c_j k u^{-(\tau-3)/2})) \nn\\
	&\le -c u^{4-\tau} k^2 \sum_{j\ge 2:c_j<\frac{1}{y_k}}  c_j^3
	\le -c k^2 \left(y_k/u\right)^{\tau-4}.
	}
We observe that $y_k/u= 8 k u^{-(\tau-1)/2}$. Consequently, $\log r_t(k) \le - c k^{\tau-2} u^{(\tau-1)(4-\tau)/2}$ which is less than  $-c|k|^{\tau-2} $. Therefore, for any $k\in \mathbb{R}$ and $u\ge 1$,
	\eqn{
	\lbeq{rt-xi-bd}
	r_t(k) \le \e^{-c|k|^{\tau-2}}.
	}
We first apply it to $t=u$. By dominated convergence, we deduce that, for $s=o(u^{\tau-3}/2)$
	\eqn{
	\lim_{u\to\infty} \int_{-\infty}^\infty \widetilde \expec[\e^{2\ii \pi k  u^{-(\tau-3)/2}\SS_u}]
	\e^{-2\ii \pi k u^{-(\tau-3)/2} s}dk
	=
	\int_{-\infty}^\infty \e^{ - (2\pi k)^2I_{\sss V}(1)/2} dk = (2\pi I_{\sss V}(1))^{-1/2}.
	}
Going back to (\ref{eq:density_inverse}) yields that
	\eqn{
	\widetilde f_{\SS_u}(s)= u^{-(\tau-3)/2}(2\pi I_{\sss V}(1))^{-1/2}(1+o(1))
	}
uniformly in $s=o(u^{\tau-3}/2)$. Furthermore, using
\eqref{eq:density_inverse}
and \refeq{rt-xi-bd}, for $t\in[u/2,u]$
	\eqn{
	\Big|\int_{-\infty}^\infty
	\widetilde \expec[\e^{2\ii \pi k  u^{-(\tau-3)/2}\SS_t}]\e^{-2\ii \pi k u^{-(\tau-3)/2} s}dk \Big|
	\le \int_{-\infty}^\infty \e^{-c|k|^{\tau-2}}dk
	}
which yields that $\widetilde f_{\SS_t}(s) \le c u^{-(\tau-3)/2}$ for all $s\in\mathbb{R}$, $u\ge 1$ and $t\in[u/2,u]$.
\qed

\section{Sample-path large deviations: proof of Theorems \ref{thm-tail-Su} and \ref{thm-sample-path-LD}}\label{sec:proofld}

\paragraph{Proof of Theorem \ref{thm-tail-Su}.} We use tilting and rewrite using $\theta=\theta_u^*$
	\eqn{
	\prob(\SS_u>0)=\phi(\theta) \widetilde{\expec}[\e^{-\theta u\SS_u} \indic{\SS_u>0}].
	}
Proposition \ref{prop-asy-phi} identifies the asymptotics of $\phi(\theta)$. What remains to do is to show that
	\eqn{
	\widetilde{\expec}[\e^{-\theta u\SS_u} \indic{\SS_u>0}]=\frac{D}{u^{(\tau-1)/2}} (1+o(1)).
	}
For this, we identify
	\eqn{
	\widetilde{\expec}[\e^{-\theta u\SS_u} \indic{\SS_u>0}]
	=u^{-1}\int_0^{\infty} \e^{-\theta v} \tilde{f}_{\sss S_u}(v/u) dv
	=u^{-(\tau-1)/2} \int_0^{\infty} \e^{-\theta v} u^{(\tau-3)/2}\tilde{f}_{\sss S_u}(v/u) dv.
	}
We use dominated convergence. For this, we use that $\theta=\thetastar_u$ converges to $\thetastar>0$ by
\eqref{thetastar-asymp} in Lemma \ref{lem-approx-VP}. Further,
by Proposition \ref{lem-dens-f}, $u^{(\tau-3)/2}\tilde{f}_{\sss S_u}(v/u)\rightarrow B$ for every $v$ fixed,
while also $u^{(\tau-3)/2}\tilde{f}_{\sss S_u}(v/u)$ is uniformly bounded. Take $u$ so large that $\thetastar_u>\thetastar/2$. Then,  $\e^{-\theta v} u^{(\tau-3)/2}\tilde{f}_{\sss S_u}(v/u)$ converges pointwise to $B\e^{-\thetastar v}$, and it is uniformly bounded by $K\e^{-\thetastar v/2}$. Then, dominated convergence yields that
	\eqn{
	\label{conv-tilde-expec}
	u^{(\tau-1)/2}\widetilde{\expec}[\e^{-\theta u\SS_u} \indic{\SS_u>0}]\rightarrow
	B\int_0^{\infty} \e^{-\thetastar v} dv=B/\thetastar.
	}
This proves Theorem \ref{thm-tail-Su} and identifies $D=B/\thetastar$.
\qed

\paragraph{Proof of Theorem \ref{thm-sample-path-LD}.} Fix $a\in[0,1]$. The case $a=0$ is obvious. We rewrite for $a \in (0,1]$ with $\theta=\theta_u^*$
	\eqn{
	\prob\big(\big|\SS_{au}-u^{\tau-2}I_{\sss E}(a)\big|>\vep u^{\tau-2}\mid\SS_u>0)
	=\frac{\widetilde{\expec}\big[\e^{-\theta u\SS_u} \indic{\SS_u>0}\indic{|\SS_{au}-u^{\tau-2}I_{\sss E}(a)|>\vep u^{\tau-2}}\big]}{\widetilde{\expec}[\e^{-\theta u\SS_u} \indic{\SS_u>0}]}.
	}
The asymptotics of the denominator were derived in \eqref{conv-tilde-expec}.
We then bound
	\eqn{
	\widetilde{\expec}\big[\e^{-\theta u\SS_u} \indic{\SS_u>0}
	\indic{|\SS_{au}-u^{\tau-2}I_{\sss E}(a)|>\vep u^{\tau-2}}\big]
	\leq 	
	\widetilde{\prob}\big(\big|\SS_{au}-u^{\tau-2}I_{\sss E}(a)\big|>\vep u^{\tau-2}\big).	
	}
By Lemma \ref{lem-expec-S}, $\widetilde{\expec}[\SS_{au}]=u^{\tau-2}I_{\sss E}(a)+o(u^{\tau-2})$, so that it suffices to prove that
	\eqn{
	\widetilde{\prob}\big(\big|\SS_{au}-\widetilde{\expec}[\SS_{au}]\big|
	>\vep u^{\tau-2}\big)=o(u^{-(\tau-1)/2}).
	}

We make crucial use of Proposition \ref{prop-laplace}, where we take $\lambda_1=\lambda, \lambda_2=0$ fixed and $t=au$, so that
	\eqn{
	\widetilde \expec\Big[\e^{ \lambda {\SS_{au} -\widetilde \expec[\SS_{au}] \over
	\sqrt{I_{\sss V}(a) u^{\tau-3}}}}\Big]
	= \e^{\lambda^2/2+ \Theta},
	}
where $|\Theta|\le o_u(1)$ since $t=au, a \in (0,1]$ fixed. By the Chernoff bound,
	\eqn{
	\widetilde{\prob}\big(\SS_{au}-\widetilde{\expec}[\SS_{au}]
	>\vep u^{\tau-2}\big)
	\leq \e^{-\vep u^{\tau-2}/\sqrt{I_{\sss V}(a) u^{\tau-3}}}
	\widetilde \expec\Big[\e^{{\SS_{au} -\widetilde \expec[\SS_{au}] \over
	\sqrt{I_{\sss V}(a) u^{\tau-3}}}}\Big].
	}
Since $\tau>3$ and since the power of $u$ appearing in the exponential equals $u^{(\tau-2)-(\tau-3)/2}=u^{(\tau-1)/2}$, this is $o(u^{-(\tau-1)/2})$.
The same proof works for $\widetilde{\prob}\big(\SS_{au}-\widetilde{\expec}[\SS_{au}]
< -\vep u^{\tau-2}\big)$ by taking $\lambda=-1$ instead.
\qed

\paragraph{Acknowledgements.}
The work of EA, RvdH, SK was supported
in part by the Netherlands Organisation for Scientific Research (NWO).
The work of JvL was supported by the European Research Council (ERC). We thank A.J.E.M. Janssen for pointing out identity \eqref{rz}.

\bibliographystyle{plain}

\end{document}